\documentclass[12pt,a4paper]{article}
\usepackage{amsmath,amstext,amssymb,amscd, color}
\usepackage{graphicx}

\usepackage[english]{babel}

\newcommand{\Xcomment}[1]{}

\makeatletter \@addtoreset{equation}{section} \makeatother

{\hfill$\qed$\medskip}

\def\qed{ \ \vrule width.1cm height.3cm depth0cm}

 \makeatletter
\renewcommand{\section}{\@startsection{section}{1}{0pt}%
{-3.5ex plus -1ex minus -.2ex}{2.3ex plus .2ex}%
{\normalfont\Large}}
 \makeatother

 \makeatletter
\renewcommand{\subsection}{\@startsection{subsection}{2}{0pt}%
{-3.0ex plus -1ex minus -.2ex}{-1.5ex plus .2ex}%
{\normalfont\large\bf}}
 \makeatother

 \Xcomment{
 \makeatletter
\renewcommand{\subsection}{\@startsection{subsection}{2}{0pt}%
{-3.0ex plus -1ex minus -.2ex}{1.5ex plus .2ex}%
{\normalfont\large\bf}}
  }
 \makeatother

\def\BA{{\mathbb{A}}}
\def\BB{{\mathbb{B}}}
\def\BC{{\mathbb{C}}}

\def\BQ{{\mathbb{Q}}}

\newcommand{\bL}{\mathbf{L}}
\newcommand{\bR}{\mathbf{R}}

\def\CA{{\cal A}}

\def\CB{{\cal B}}

\def\CC{{\cal C}}

\def\CV{{\cal V}}

\def\ba{{\bf a}}

\def\bbf{{\bf f}}
\def\bk{{\bf k}}

\def\tilde{\widetilde}
\def\hat{\widehat}

\def\lra{\longrightarrow}

\begin{document}


\title{Quaternities, correspondences, and tetrahedron equations (Summa tetralogiae)}


\author{Gleb Koshevoy, Vadim Schechtman, and Alexander Varchenko}

\date{\today}

\maketitle

{\em La Nature est un temple o\`u de vivants piliers 

Laissent parfois sortir de confuses paroles; 

L'homme y passe \`a travers des for\^ets de symboles 

Qui l'observent avec des regards familiers.}

\

{\em Ch. Baudelaire, Correspondances}\footnote{Note the resemblance of the structure of this sonnet to the structure of our tetrahedron equation, cf. (\ref{tetr-eq-corresp}).}

\section{Introduction}

\subsection{•} The aim of this note is:

\

(a) to propose a generalization of tetrahedron equations from \cite{S} and of their solutions. Due to appearance of a larger number of parameters the  $R$-matrices from \cite{S}  will be replaced by "$R$-correspondences".

\

(b) To rephrase these equations in terms of Wronskian evolutions in the spirit  
of \cite{SV}.

\

(c) To discuss some elementary structures of cohomological flavour   lying behind our 
considerations. We call them "quaternities", or "bibitorsors"; they  might 
be not without an independent interest.

\subsection{Sonnet equations.} They are a refined version of tetrahedron equations.  To describe them we need some combinatorics 
which is a particular case of \cite{MS}. 

{\it Categories $\CA(3,2)$ and $\CA(4,2)$.} Consider the symmetric group $S_3$  
with Coxeter generators $s_1, s_2$. The relations are 
$$
s_1s_2s_1 = s_2s_1s_2
$$
and $s_1^2 = s_2^2 = e$. The set $A(3,2)$ of reduced decompositions of the longest element 
in $S_3$ consists of two elements $s_1s_2s_1, s_2s_1s_2$, to be denoted $(121), (212)$.      

The category $\CA(3,2)$ has two objects $(121), (212)$ and one arrow 
$R:\ (121) \lra (212)$.

\

Consider the symmetric group $S_4$ with Coxeter generators $s_1, s_2, s_3$. The relations are 
$$
s_is_{i+1}s_i = s_{i+1}s_is_{i+1}, \ i = 1, 2,
$$
and $s_i^2 = e, \ 1\leq i\leq 2$. As before, products $s_{i_1}\ldots s_{i_k}$ will be denoted $(i_1\ldots i_k)$. Let $A(4,2)$ denote the set of all reduced decompositions 
of the longest element in $S_4$. For example 
$$
(121321), (321323) \in A(4,2). 
$$
By definition $Ob \CA(4,2) = A(4,2)$. 
Altogether $A(4,2)$ consists of $14$ elements; let us describe them together with 
the morphisms in $\CA(4,2)$. By definition morphisms in $\CA(4,2)$ are compositions of 
elementary arrows. Elementary arrows are of two types:

\

(a) {\em $R$-type}: these are arrows of the form
$$
R(k):\ (... i_k i_{k+1} i_{k+2} ...) \lra (... i'_k i'_{k+1} i'_{k+2} ...), \ 
1\leq k\leq 4,
$$
where $(i_k, i_{k+1}, i_{k+2}) = (j, j+1, j)$ and 
$(i'_k, i'_{k+1}, i'_{k+2}) = (j+1, j, j+1)$ for $j = 1$ or $2$, the other elements 
not being changed. 

For example
$$
R(1): \ (121321) \lra (212321).
$$

(b) {\em $L$-type}: these are arrows of the form
$$
L(k):\ (... i_k i_{k+1} ...) \lra (... i_{k+1} i_k ...),\ 1\leq k \leq 5,
$$
identical on other elements, where $(i_k, i_{k+1}) = (1, 3)$ or $(3, 1)$. 

For example
$$
L(3):\ (121321) \lra (123121). 
$$
Now we can draw all objects in $\CA(4,2)$, and almost all arrows. 

Namely, there are two chains of elementary arrows: the first (upper) one is
$$
(121321) \overset{R(1)}\lra (212321) \overset{R(3)} \lra (213231)  
\overset{L(2)}\lra (231231) \overset{L(5)}\lra (231213) 
$$
$$
\overset{R(3)}\lra (232123) \overset{R(1)}\lra (323123)  
\overset{L(3)}\lra (321323)
$$
Let $C_+$ denote the composition 
$$
C_+ = L(3)R(1)R(3)L(5)L(2)R(3)R(1):\ (121321) \lra (321323).
$$
The second (lower) chain is 
$$  
(121321) \overset{L(3)}\lra (123121) \overset{R(4)} \lra (123212)  
\overset{R(2)}\lra (132312) \overset{L(4)}\lra (132132) 
$$
$$
\overset{L(1)}\lra (312132) \overset{R(2)}\lra (321232)  
\overset{R(4)}\lra (321323)
$$
Let $C_-$ denote the composition 
$$
C_- = R(4)R(2)L(1)L(4)R(2)R(4)L(3):\ (121321) \lra (321323).
$$
We impose the relations
$$
L(i)^2 = \text{Id},\ 1\leq i \leq 5,
$$
$$
L(i)L(j) = L(j)L(i)\ \text{if}\ |i - j| \geq 2,
$$
and 
$$
C_+ = C_-
$$
(the sonnet relation), which may be written symbolically in the form
$$
LRRLLRR = RRLLRRL.
$$
If we forget the letter $L$ we get an equation of the form
$$
RRRR = RRRR,
$$
the tetrahedron relation.

\

Let $\CC$ be a category. Let us call a {\it sonnet equation in $\CC$} a functor 
$\BA:\ \CA(4, 2) \lra \CC$. So it is a collection of $14$ objects $\BA(\ba)\in Ob \CC, \ba\in A(4,2)$ and $14$ morphisms $\BA(f):\ \BA(\ba) \lra \BA(\ba')$ given for each elementary arrow 
$f:\ \ba \lra \ba'$, i.e. $f = R(i)$ or $f = L(j)$, such that $\BA(L(j))^2 = \text{Id}$, 
$\BA(L(j))\BA(L(j')) = \BA(L(j'))\BA(L(j))$ for $|j - j'| > 1$, and
$$
L(3)R(1)R(3)L(5)L(2)R(3)R(1) = R(4)R(2)L(1)L(4)R(2)R(4)L(3):\ 
$$
$$
\BA(121321) \lra \BA(321323).
$$
where we have denoted for brevity $\BA(f)$ by $f$. 

\

More generally, for any $n\geq 2$ we can define in a similar manner the sets $A(n,2)$ 
of all reduced decompositions of the longest element in $S_n$, and a category $\CA(n,2)$, 
cf. \cite{MS}. 

\ 

\subsection{Category of rational correspondences.}
In the examples below the target category $\CC$ will be the category whose objects are   
varieties over a fixed base field, or, more generally, schemes, and morphisms are rational correspondences. 

Namely, let $X$, $Y$ be two varieties. 
Let us call a rational map $f: X\lra Y$ a map of varieties $f_0: U \lra Y$ where 
$U\subset X$ 
is a dense Zarisky open in $X$, where two maps $f_0:\ U \lra Y$ and $f_0':\ U' \lra Y$ are identified if they coincide on $U\cap U'$. 

Let us say that $f$ is {\it  good} if  
$f_0(U)$ contains a Zarisky dense open in $Y$. Good rational maps may be composed.

A rational correspondence between $X_1$ and $X_2$ is a rational map
$$
P \lra X_1\times X_2
$$
whose components $p_i: P \lra X_i,\ i = 1, 2$ are good. 

Rational correspondences may be composed in the usual way. Namely, given another correspondence 
$$
(q_2, q_3):\ Q \lra X_2 \times X_3,
$$
let $R = P\times_{X_2} Q$; let $r_1: R\lra P$, $r_3: R\lra Q$ be the projections.
The composition $QP: R \lra X_1 \times X_3$ is defined as $(p_1r_1, q_3r_3)$.

The composition is associative, and $\text{Id}_X$ is the diagonal 
$X \lra X\times X$. This way we get the category $\CV$ of varieties and good rational 
correspondences. It will be our target category.  

\subsection{Categories $\CB(n,2)$.} 

\

More generally, for each $n\geq 2$ our examples give rise to functors 
$\CA(n,2) \lra \CV$.
In \cite{MS} the sets $B(n,2)$ are introduced as certain quotients of $A(n,2)$. For example $B(4,2)$ is obtained from $A(4,2)$ by identifying elements related by $L$ transformations. It contains $8$ elements. We can define the category 
$\CB(n,2)$ with $Ob \CB(n,2) = B(n,2)$ similarly to $\CA(n,2)$.

\ 

A {\it tetrahedron equation} with values in a category $\CC$ is a functor
$$
\CB(4,2) \lra \CC
$$
It will have a form 
$$
RRRR = RRRR
$$
We have the canonical projection functor
$$
\CA(4,2) \lra \CB(4,2),
$$
and one sees immediately that it is an equivalence of categories. Therefore 
 tetrahedron equations are equivalent to sonnet equations. 
 
\subsection{•} Our examples are given in Sections \ref{sect-smaller} 
and \ref{flacon}. 
The spaces $\BA(\ba)$ therein are some affine spaces $\BA^N$.

\subsection{Evolutions along reduced decompositions.} Fix $n\geq 2$. For a matrix 
$$
A = \left(\begin{matrix}
a & b\\
c & d
\end{matrix}\right)
$$
and $1 \leq i \leq n - 1$ let $\phi_i(A)$ denote the $n\times n$ matrix of the linear 
transformation

\

$$
\phi_i(A)(e_i) = ae_i + be_{i+1},\ \phi_i(A)(e_{i+1}) = ce_i + de_{i+1}, 
$$
$$ 
\phi_i(A)(e_j) = e_j\ \text{for}\ j\neq i, i+1.
$$
Pick some Coxeter generators $s_1, \ldots, s_{n-1}$ of $S_n$; let $w_0\in S_n$ be the longest element. To a reduced decomposition
\begin{equation}\label{reduced-decomp}
w_0 = s_{i_1}\ldots s_{i_{c(n)}},\ c(n) = n(n-1)/2, 
\end{equation} 
and a collection of matrices $A_k\in Mat_2$ let us assign matrices
$$
B_k = \phi_{i_1}(A_1)\ldots \phi_{i_k}(A_k)\in Mat_n, \ 1\leq k \leq c(n).
$$
We will call the sequence $B_1, \ldots, B_{c(n)}$ {\em the evolution along the decompostion} (\ref{reduced-decomp}).

\

\subsection{$c$-upper triangular version.}\label{c-upper-evol} Let us call a square matrix 
$M = (a_{ij})$ {\it $c$-upper triangular} (or standard $c$-Borel) if $a_{ij} = 0$ for $i < j$. Let us denote 
by $B^c_n \subset Mat_n$ the subspace of $c$-upper triangular matrices.

\

It turns out that for every reduced decomposition as above, if all $A_i\in B^c_2$ 
then all $B_i\in B^c_n$.

\

\subsection{Wronskian evolutions.}
On the other hand we can realize the standard $n$-dimentional affine space as the space of 
polynomials $\bk[x]_{\leq n-1}$ of degree $\leq n - 1$ over the ground field $\bk$ 
which we suppose here to be of characteristic $0$. So an $n\times n$ matrix $M$ is 
the same as a sequence 
$$
\bbf = \bbf(M) = (f_1(x), \ldots, f_n(x))
$$
of such polynomials. Let us assign to $\bbf$ a sequence of Wronskians
$$
Wr(\bbf) := (Wr(f_1), Wr(f_1, f_2), \ldots, Wr(f_1, \ldots, f_n))
$$

\

Let $B_i, 1\leq i \leq c(n)$ be as in \ref{c-upper-evol}, whence the  sequence 
$\bbf(B_1), \ldots, \bbf(B_{c(n)})$. 

In Section \ref{sec-wronskian-ev} we describe the corresponding sequence of Wronskians 
$$
Wr(\bbf(B_1)), \ldots, Wr(\bbf(B_{c(n)})).
$$
For the usual upper triangular case this has been done in \cite{SV}.

\subsection{•} The formalism of quaternities/bitorsors is described in Section 
\ref{sect-quater}.

\subsection{•} Correspondences related to tetrahedron equations have been discussed in 
\cite{IKR}. 

\subsection{Acknowledgements.} G.K. thanks IHES for the hospitality and perfect working conditions.  
A.V. was supported in part by Simons Foundation grant TSM-00012774. 

\section{Berenstein - Zelevinsky transitions}

\subsection{Flags and bases.} Let $V = V_{r+1} = \BC[x]_{\leq r}$ denote the space of complex polynomials of degree $\leq r$. Let $Fl_r$ denote the variety of flags of linear subspaces 
$$
L_\bullet:\ L_1\subset L_2\subset\ldots \subset L_r\subset L_{r+1} = V_{r+1},\ \dim L_i = i.
$$
Let   $U\subset Fl_r$ be the Zarisky dense open subvariety 
of flags in general position. The meaning of the words "in general position" will be 
as follows:

\

(1) if $u_1 = u_1(x)\in L_1$ is a base vector then we demand that $u_1(0)\neq 0$;

(2) pick $u_1(x)$ as above; we demand that there exists $u_2(x)\in L_2$ such that 
$u_1, u_2$ form a base of $L_2$, $u_2(0) = 0$, and $u_2'(0)\neq 0$; 

(3) pick $u_1, u_2$ as above; we demand that there exists $u_3(x)\in L_3$ such that 
$u_1, u_2, u_3$ form a base of $L_3$, $u_3(0) = u'_3(0) = 0$, and $u''_3(0)\neq 0$;

etc.

We will also call flags in general position {\it admissible flags}. 

\

Let $L_\bullet\in U$. There exists a unique collection $u_\bullet = u_\bullet(L_\bullet) 
= \{ u_1, \ldots, u_r\} \in V$ such that 

\

(a) $u_1, \ldots, u_i$ form a base of $L_i$;

\

(b) $u_1(0) = 1, u_2 = x +\ldots, u_3 = x^2/2 + \ldots$, etc., i.e. 
$u_i = x^{i-1}/(i-1)! + \ldots$ where $\ldots$ mean the terms of higher order.

\

We will call $u_\bullet$ the {\it distingushed base} of $L_\bullet$.

\subsection{Adding up the diagonal.} More generally, given a collection $c = \{c_1, \ldots, c_{r+1}\in \BC^*\}$, there exists a unique collection $u_c = \{u_1, \ldots, u_{r+1}\}$ such that 

\

(a) $u_1, \ldots, u_i$ form a base of $L_i$;

\

$(b)_c$ $u_i = c_ix^{i-1}/(i-1)! + \ldots$.

\

Note that $u_{r+1} = c_rx^r/r!$.

\subsection{Normalized collections and transformations $A_i(a)$.} Let us call a {\em normalized collection} with parameters $c_1, \ldots, c_{r+1},\ c_i\in \BC^*$ 
a collection $\{ u_1, \ldots, u_{r+1}\}\subset V$ such that 
$$
u_i = c_ix^{i-1}/(i-1)! + \ldots,\ i = 1, \ldots, r+1,
$$
where $...$ denotes the terms of higher order. 

\

We may identify $V^{r+1}$ with the space $M = Mat_{r+1}(\BC)$ of $(r+1)\times (r+1)$ matrices; 
its elements may be written as columns
$$
\left(\begin{matrix}
u_1 \\ . \\ . \\ . \\ u_{r+1}
\end{matrix}\right)
$$

\

Let $U\subset V^{r+1}$ denote the variety of 
all normalized collections. It may be identified  with the standard upper Borel subgroup of $GL_{r+1}(\BC)$.

\

For $a\in \BC^*$ let
$$
A(a) = \left(\begin{matrix}
a^{-1} & 1\\
0 & a
\end{matrix}\right)
$$
Let us define operators $A_i(a):\ M\lra M$ by 
$$
A_i(a)\left(\begin{matrix}
u_1\\ \ldots \\ u_i\\ u_{i+1}\\ \ldots\\ u_{r+1}
\end{matrix}\right)  = \left(\begin{matrix}
u_1\\ \ldots \\ a^{-1}u_i + u_{i+1}\\ au_{i+1}\\ \ldots\\ u_{r+1}
\end{matrix}\right),\   
$$
$i = 1, \ldots, r$, cf. \cite{BZ} (4.9).  

\

We see that $A_i(U)\subset U$.

\subsection{Theorem}\label{bz-map} {\em (braid relation)}, cf. \cite{BZ}, (5.2). {\em We have
$$
A_1(c)A_2(b)A_1(a) = A_2(c')A_1(b')A_2(a')
$$
where
$$
a' = \frac{bc}{ac+b},\ b' = ac,\ c' = \frac{ac+b}{c}.
$$} 

$\Box$

\section{Wronskians and transition maps}

\subsection{Wronskians.} For a sequence $u = (u_1(x), \ldots, u_i(x)),\ u_j(x)\in \BC[x]$, define the Wronskian matrix
$$
\mathcal{W}r(u_1, \ldots, u_i)  = \left(\begin{matrix}
u_1 & u_2 & \ldots & u_i\\
u'_1 & u'_2 & \ldots & u'_i\\
. & . & \ldots & .\\
u^{(i-1)}_1 & u^{(i-1)}_2 & \ldots & u^{(i-1)}_i\\
\end{matrix}\right),
$$
and let $Wr(u) = Wr(u_1,\ldots, u_i)$ denote its determinant.

\

For example 
$$
\mathcal{W}r(1, x, x^2/2, \ldots, x^{i-1}/(i-1)!)
$$
is an upper triangular matrix with $1$'s on the diagonal,  
whence 
$$
Wr(1, x, x^2/2, \ldots, x^{i-1}/(i-1)!) = 1.
$$
On the other hand
$$
\mathcal{W}r\left(1, x, \frac{x^2}{2}, \ldots , \frac{x^{i-1}}{(i-1)!}, 
\frac{x^{i+1}}{(i+1)!}\right) 
$$
is an upper triangular matrix whose diagonal part is $Diag(1, \ldots, 1, x)$,  
whence
$$
Wr\left(1, x, \frac{x^2}{2}, \ldots , \frac{x^{i-1}}{(i-1)!}, 
\frac{x^{i+1}}{(i+1)!}\right) = x.
$$

\subsection{Wronskian of a normalized collection.} Let $u = (u_1, \ldots, u_{r+1})$ be a normalized collection with 
parameters $\{ c_i\}$. Denote
$$
f_i(x) := Wr(u_1,\ldots, u_i)
$$
and 
$$
Wr(u) := (f_1, \ldots, f_{r+1}).
$$
We have 
$$
f_i(x) = c_1\cdot\ldots\cdot c_i + O(x).
$$
Denote 
$$
a_i := \prod_{j=1}^i c_j.
$$
We have
$$
\frac{a_{i-1}a_{i+1}}{a_i} = \left(\prod_{j=1}^{i-1}c_j\right)c_{i+1}.
$$

\subsection{Proposition.} {\em We have 
$$
Wr(A_i(a)u) = (f_1,\ldots, f_{i-1}, a^{-1}f_i + \frac{a_{i-1}a_{i+1}}{a_i} x + O(x^2), 
af_{i+1}, f_{i+2}, \ldots, f_{r+1}).
$$}

\subsection{Wronskian evolution} We can go in the opposite direction. Let us consider a 
variety   
$$
W = W_{r+1} = \{
f = (f_1, \ldots, f_{r+1}),\ f_j\in \BC[x],   f_j(0)\neq 0,\ 1\leq j\leq r+1\}. 
$$
We have defined the Wronskian map
$$
Wr:\ U\lra W.
$$
Let $f\in W,\ f_j(0) = a_j$.
There exists a unique solution $\hat{f}_i$ of the differential equation 
$$
Wr(f_i, \hat{f}_i) = f_{i-1}f_{i+1}
$$
of the form
$$
\hat{f}_i = \frac{a_{i-1}a_{i+1}}{a_i} x + O(x^2).
$$
For $a\in \BC^*$ denote 
$$
\tilde{f}_i  = a^{-1}f_i + \hat{f}_i 
$$
and define
$$
A_i(a)f = (f_1,\ldots,f_{i-1},\tilde{f}_i, af_{i+1}, f_{i+2},\ldots, f_{t+1}).
$$

\subsection{Theorem} {\em For $u\in U$ 
$$
Wr(A_i(a)u) = A_i(a)Wr(u).
$$}

\

Cf. \cite{SV}, Thm. 4.4, (4.8.1).

\subsection{Examples} Let us start with 
$$
u = (1, x, \frac{x^2}{2},\ldots, \frac{x^r}{r!}), 
$$
so 
$$
f := Wr(u) = (1, 1, \ldots, 1).
$$
Then 
$$
A_2(a)u = (1, x, a^{-1}x + \frac{x^2}{2}, \frac{ax^2}{2}, \frac{x^3}{3!},\ldots, \frac{x^r}{r!}),
$$ 
and
$$
A_2(a)f = (1, a^{-1} + x, a, 1, \ldots, 1).
$$
Next, 
$$
A_1(b)A_2(a)u = (b^{-1} + a^{-1}x + \frac{x^2}{2}, b(a^{-1}x + \frac{x^2}{2}), 
\frac{ax^2}{2}, \frac{x^3}{3!},\ldots, \frac{x^r}{r!}), 
$$
and 
$$
A_1(b)A_2(a)f = (b^{-1} + a^{-1}x + \frac{x^2}{2}, b(a^{-1} + x), a, 1, \ldots, 1).
$$

\section{Sergeev $(\alpha)$ as a Wronskian evolution}.

\subsection{Identities in $SL_n(\BC)$.} For $1\leq i,j\leq n$, $i\neq j$, let us denote by $e^a_{ij}\in SL_n(\BC)$ the matrix $I_n + ae'_{ij}$ where 
$e'_{ij}$ is the matrix with $1$ in the $i$-th row and the $j$-th column and $0$ everywhere else, cf \cite{M}, \S 5. We have
$$
(e^a_{ij})^{-1} = e^{-a}_{ij},
$$
and
$$
[e^a_{ij}, e^b_{jk}] := e^a_{ij}e^b_{jk}e^{-a}_{ij}e^{-b}_{jk} = e^{ab}_{ik}  
$$
if $i\neq k$.
For example  
$$
[e^a_{12}, e^b_{23}] = e^{ab}_{13}.  
$$

\

These identities allow us to define the Wronskian evolution.

\subsection{A triple evolution.}\label{triple} Denote
$$
A(a) = \left(\begin{matrix}
a & -1\\
1 & 0
\end{matrix}\right),
$$
and let $A_i(a)$ denote the $(n\times n)$ block matrix acting as $A(a)$ on the $i$-th and $(i+1)$-th row and column, and as identity elsewhere.
We have
$$
A_1(a)A_2(b) = \left(\begin{matrix}
a & -b & 1\\
1 & 0 & 0\\
0 & 1 & 0
\end{matrix}\right). 
$$
Note that if we interchange the last two rows of the above matrix we get a matrix which is upper triangular with respect to the complementary diagonal.  

We have
$$
A_1(a)A_2(b)A_1(c) = \left(\begin{matrix}
ac - b & -a & 1\\
c & -1 & 0\\
1 & 0 & 0
\end{matrix}\right). 
$$

\subsection{Notation.} Let $I(m)$ denote the $(m\times m)$-matrix with elements $d_{ij} = 1$ if $i+j=m+1$ and $0$ otherwise. 

\subsection{Triple evolution (cont-d).}\label{triple-cont} Let $u_i(x)\in \BC[x]$, $i = 1, 2, 3$, 
$$
u = \left(\begin{matrix}
u_3 \\ u_2 \\ u_1\\
\end{matrix}\right) = \left(\begin{matrix}
x^2/2 + \ldots \\ x + \ldots \\ 1 + \ldots\\
\end{matrix}\right). 
$$
We have 
$$
u = I(3)\left(\begin{matrix}
u_1 \\ u_2 \\ u_3\\
\end{matrix}\right).
$$ 
We have
$$
A_1(a)A_2(b)A_1(c)u = \left(\begin{matrix} (ac-b)u_3 - au_2 + u_1\\
cu_3 - u_2 \\ 
u_3
\end{matrix}\right) = 
$$
$$
=  \left(\begin{matrix}  1 + \ldots \\ - x + \ldots \\ x^2/2 + \ldots\\
\end{matrix}\right).
$$
This evolution corresponds to the decomposition
$$
w_0 = s_1s_2s_1
$$
of the longest element of $S(3)$.

\subsection{A quadruple evolution.}\label{quadruple} Similarly consider the decomposition
$$
w_0 = s_1(s_2s_1)(s_3s_2s_1)
$$
of the longest element in $S(4)$. The corresponding product in $SL(4)$ will be
$$
A(a,b,c) := A_1(a)A_2(b)A_1(a)A_3(c)A_2(b)A_1(a).
$$
It is an upper triangular matrix with respect to the complementary diagonal with $\pm 1$'s 
on the complementary diagonal. So
\begin{equation}\label{quadruple-eq}
A(a,b,c)\left(\begin{matrix}
u_4\\ u_3\\ u_2\\ u_1
\end{matrix}\right) = 
\left(\begin{matrix}
\pm 1 + \ldots\\
\pm x + \ldots\\
\pm x^2/2 + \ldots\\
\pm x^3/6 + \ldots
\end{matrix}\right).
\end{equation}

\section{Quaternity calculus}\label{sect-quater}

\

{\em To be or not to be?}

\

{\em Hamlet}

\subsection{Some terminology and notations.} The ground ring $\BC$ will be the field of 
complex numbers, or more generally any commutative ring. In the Wronskian Section we will suppose that it contains $\BQ$.

The space of $(n\times n)$-matrices will be denoted $Mat_n = Mat_n(\BC)$; $I_n\in Mat_n$ will denote the unity matrix.  

\

A square matrix $A = (a_{ij})$ will be called upper (resp. lower) triangular if $a_{ij} = 0$ for $i > j$ (resp. for $i < j$). The group of upper 
(resp. lower) triangular $(n\times n)$-matrices with invertible elements on the diagonal  will be denoted 
$B_+(n)$ (resp. $B_-(n)$). The subgroups of matrices with $1$'s on the diagonal will be denoted by $N_+(n), N_-(n)$.

\

Let $A = (a_{ij})$ be an $(n\times n)$-matrix. It {\it complementary diagonal} is the set 
of elements $a_{ij}$ with $i+j=n+1$. We will denote by $B^c_+(n)$ (resp. by $B^c_-(n)$) 
the set of matrices with nonzero elements only above (resp. under)  or on the complementary  diagonal, with invertible invertible elements on the complementary diagonal. 
Elements of $B^c_+(n)$ (resp. of $B^c_-(n)$) will be called $c$-upper- (resp. $c$-lower-) 
triangular matrices.  

We will denote by $N^c_\pm(n)\subset B^c_\pm(n)$ the subsets of matrices with $1$'s on the 
complementary diagonal; $I(n)$ will denote the matrix with $1$'s on the complementary diagonal and zeros elsewhere.  

Thus $N^c_+(n)\cap N^c_-(n) = \{ I(n)\}$. 

\

We have an injective group homomorphism
$$
\iota: S(n)\lra Mat_n
$$
sending $w\in S(n)$ to the corresponding permutation matrix. 

Thus 
$I(n) = \iota(w_0(n))$ where $w_0(n)$ is the longest element of $S(n)$.

$e(n)$ will denote the unit element of $S(n)$. 

\

If $A\in Mat_n$ then we will denote by $L_A: Mat_n\lra Mat_n$ (resp. by 
$R_A: Mat_n\lra Mat_n$) the operator $L_A(B) = AB$ (resp. $R_A(B) = BA$).

\subsection{Some matrix identities.}

(1) $w_0(n)^2 = e(n)$, therefore $I(n)^2 = I_n$, therefore 
$L_{I(n)}^2 = R_{I(n)^2} = \text{Id}_{Mat_n}$. 

\

(2) If $A\in Mat_n$ then $I(n)A$ (resp $AI(n)$) 
is $A$ with the order of lines (resp. of columns) reversed. 

\subsection{Quaternity.}
The operators $L = L_{I(n)}, R = R_{I(n)}$ induce four bijections which fit in a square called  {\em $B$-quaternity}, cf. \cite{J}\footnote{Quaternity:
An image with a four-fold structure, usually square or circular and symmetrical; psychologically, it points to the idea of wholeness (Carl Jung).}
\begin{equation}\label{quat-diagr}
\begin{matrix}
B_+ & \overset{L}{\lra} & B_-^c\\ 
R\uparrow & & \downarrow R\\
B_+^c & \overset{L}{\longleftarrow} & B_-
\end{matrix}
\end{equation}
The composition of all four, starting from any place, is identity.

\

Similarly we define $N$-quaternity with the letter $B$ replaced by $N$. 

\subsection{Actions.} The above operators define an action of the group $B_+$ (resp. $B_-$) upon the set $B_-^c$ from the left (resp. from the right). Namely for $x\in B_+, y\in B_-^c, 
z\in B_-$ we set 
$$
xy := L(xL^{-1}(y)), \ yz := R(R^{-1}(y)z).
$$
These two actions commute, and $B_-^c$ is a left $B_+$-torsor and a right $B_-$-torsor.
 
In other words $B_-^c$ is a {\it $(B_+, B_-)$-bitorsor}, in the sense of \cite{Gi}, D\'efinition 1.5.3. 


\

Similarly $B_+^c$ is a $(B_-, B_+)$-bitorsor.

\

The whole quaternity diagram (\ref{quat-diagr}) might be called a {\em bibi-torsor}.  

\

Similarly with $B$ replaced by $N$.

\section{$R$-correspondences and tetrahedron equations ($c$-upper triangular case)}

\subsection{•} For $A\in Mat_2$ and $1\leq i\leq n-1$ we denote by $A_i\in Mat_n$ 
the block matrix acting as $A$ in the $i$-th and $(i+1)$-th rows and columns and as identity 
elsewhere.

For example if $n = 3$, 
$$
A_2 = \left(\begin{matrix} 1 & 0 & 0\\
0 & a_{11} & a_{12}\\
0 & a_{21} & a_{22}
\end{matrix}\right).
$$  

\subsection{Double and triple products.}\label{double-triple} 
Let us denote
$$
A(a,b,c) = \left(\begin{matrix}
a & b\\
c & 0
\end{matrix}\right).
$$
We have
\begin{equation}\label{double-one}
A(a_1,b_1,c_1)_1 A(a_2,b_2,c_2)_2 = \left(\begin{matrix}
a_1 & b_1a_2 & b_1b_2\\
c_1 & 0 & 0 \\
0 & c_2 & 0
\end{matrix}\right).
\end{equation}
Note that if in the last matrix we interchange the second and the third row, or the first and the second column, we get an element of $B_+^c(3)$.

Next we have 
$$
A(a_1,b_1,c_1)_1 A(a_2,b_2,c_2)_2 A(a_3,b_3,c_3)_1 = 
\left(\begin{matrix}
a_1a_3 + b_1a_2c_3 & a_1b_3 & b_1b_2\\
c_1a_3 & c_1b_3 & 0\\
c_2c_3 & 0 & 0
\end{matrix}\right).
$$
On the other hand
\begin{equation}\label{double-two}
A(a_1,b_1,c_1)_2 A(a_2,b_2,c_2)_1 = \left(\begin{matrix}
a_2 & b_2 & 0\\
a_1c_2 & 0 & b_1\\
c_1c_2 & 0 & 0
\end{matrix}\right)
\end{equation}
which becomes $c$-upper triangular after the interchange of the second and the third column, 
or of the first and the second row.

\

Next
$$
A(a_1,b_1,c_1)_2 A(a_2,b_2,c_2)_1 A(a_3,b_3,c_3)_2 = 
\left(\begin{matrix}
a_2 & b_2a_3 & b_2b_3\\
a_1c_2 & b_1c_3 & 0\\
c_1c_2 & 0 & 0
\end{matrix}\right).
$$

\subsection{$R$-correspondence.} The $R$-correspondence is a  subvariety   
$$
R\subset\BA^9 \times \BA^9.
$$
If we denote the coordinates in $\BA^9 \times \BA^9$ by
$(a_1, \ldots, c_3; a'_1, \ldots, c'_3)$ then $R$ is given by the equation
\begin{equation}\label{R-corresp}
A(a_1,b_1,c_1)_1 A(a_2,b_2,c_2)_2 A(a_3,b_3,c_3)_1 = 
A(a'_1,b'_1,c'_1)_2 A(a'_2,b'_2,c'_2)_1 A(a'_3,b'_3,c'_3)_2.
\end{equation}
Explicitely this equation is equivalent to six equations
\begin{equation}\label{RC1}
a_1a_3 + a_2b_1c_3 = a'_2
\end{equation}
\begin{equation}\label{RC2}
a_1b_3 = b'_2a'_3
\end{equation}
\begin{equation}\label{RC3}
b_1b_2 = b'_2b'_3
\end{equation}
\begin{equation}\label{RC4}
c_1a_3 = a'_1c'_2
\end{equation}
\begin{equation}\label{RC5}
c_2c_3 = c'_1c'_2
\end{equation}
\begin{equation}\label{RC6}
c_1b_3 = b'_1c'_3
\end{equation}

\subsection{Parametrization of $R$.}\label{param} Let us forget $a_2, a'_2$. The remaining $16$ variables
$$
(a_1, b_1, c_1, b_2, c_2, a_3, b_3, c_3; a'_1, b'_1, c'_1, b'_2, c'_2, a'_3, b'_3, c'_3)
$$
are related by $5$ equations (\ref{RC2}) - (\ref{RC6}). Among these equations there are 
two pairs: (\ref{RC2}) - (\ref{RC3}) and (\ref{RC4}) - (\ref{RC5}). 

If we choose in each pair $(a'_1, c'_1)$, $(a'_3, b'_3)$, $(b'_1, c'_3)$ one variable, for example $a'_1, a'_3, b'_1$ then the other primed variables are defined uniquely 
(in general position, which means on a Zarisky dense subvariety) from our equations. 

\

For example if we choose $a'_1, a'_3, b'_1$ as free parameters then 
\begin{equation}
c'_1 = c_2c_3/(a_3c_1)\cdot a'_1,\ b'_2 = a_1b_3/a'_3,\ c'_2 = a_3c_1/a'_1,  
\end{equation}
\begin{equation}
b'_3 = b_1b_2/(a_1b_3)\cdot a'_3,\ c'_3 = b_3c_1/b'_1.
\end{equation}
Thus $R$ is a $12$-dimensional rational (i.e. birationally isomorphic to $\BA^{12}$)   subvariety of $\BA^9\times \BA^9$ ("a cell"). 


\subsection{Compositions of correspondences.} We denote a correspondence 
$R\subset X\times Y$ by 
$$
X\overset{R}{\lra} Y.
$$ 
Let $R\subset X\times Y$, $S\subset Y\times Z$ 
be two correspondences. Their composition $SR\subset X\times Z$ is defined by 
$$
SR = p_{13}(R\times_Y S)
$$
where $R\times_Y S \subset X\times Y\times Z$, and $p_{13}: X\times Y\times Z\lra X\times Z$ 
is the projection. In other words
$$
SR = \{ (x,z)\in X\times Z|\ \exists y\in Y \text{such that} (x,y)\in R, (y,z)\in S\},
$$
cf. \cite{CG}, 2.7.5. 

The compoosition is associative.

This way we get a category, to be denoted by $Corr$,  whose objects are varieties and morphisms are correspondences. 
The identity Id$_X$ is the diagonal $\Delta\subset X\times X$. 

\subsection{Tetrahedron equations.} They are associated with reduced decompositions 
of $w_0(4)\in S(4)$. The set of such decompositions consists of $14$ elements; it is the set $A(4, 2)$ related to the second Bruhat order from \cite{MS}.
They can be found in \cite{SV}, (8.10).

\ 

The equation is an equality of two correspondences in $\BA^{18}\times \BA^{18}$.
There are $8$ flips of type $R$, and $6$ flips of type $L$ as in \cite{SV}, (8.14), 
and the equation has the same form as in \cite{SV}, 8.2.3, (8.16). 

\

\subsection{Problem: do both ways give the same result?} Let us see if the diagram 

\begin{equation}\label{tetr-eq-diagram}
\begin{matrix} 
 \BA(212321) & \overset{\bR(3)}\lra & \BA(213231) \overset{\bL(2)}\lra \BA(231231) 
 \overset{\bL(5)}\lra  \BA(231213) & 
 \overset{\bR(3)}\lra & \BA(232123) \\ 
\bR(1)\uparrow & & & &  \downarrow\bR(1)\\
\BA(121321) & & & & \BA(323123)\\
\bL(3)\downarrow & & & & \downarrow\bL(3) \\
\BA(123121) & & & & \BA(321323)\\
\bR(4)\downarrow & & & & \uparrow\bR(4) \\
\BA(123212) & \overset{\bR(2)}\lra & \BA(132312) \overset{\bL(4)}\lra \BA(132132) 
\overset{\bL(1)}\lra  \BA(312132) & 
\overset{\bR(2)}\lra &  \BA(321232)   
\end{matrix}
\end{equation}
in $Corr$ commutes birationally. 

\

In other words, is it true that 
\begin{equation}\label{tetr-eq-corresp}
\bL(3)\bR(1)\bR(3)\bL(5)\bL(2)\bR(3)\bR(1) \overset{\sim}{=} \bR(4)\bR(2)\bL(1)\bL(4)\bR(2)\bR(4)\bL(3)
\end{equation} 
where $\overset{\sim}{=}$ means the equality of Zarisky dense subspaces?

\

Here symbols $\BA(k_1k_2k_3k_4k_5k_6)$ denote copies of $\BA^{18}$. 

$\bL(i)$ is the graph of the function
$$
(a_i,b_i,c_i, a_{i+1},b_{i+1},c_{i+1})\mapsto (a_{i+1},b_{i+1},c_{i+1}, a_i,b_i,c_i)
$$
and identical on other coordinates.

\

In  the rest of this Section we will compute the left and the right hand sides of 
(\ref{tetr-eq-corresp}). 

\subsection{Notations.} Denote the spaces in the upper way by   
$$
\BA_1 = \BA(121321), \BA_2 = \BA(212321), \BA_3 = \BA(213231), \BA_4 = \BA(231231),
$$
$$ 
\BA_5 = \BA(231213), \BA_6 = \BA(232123), \BA_7 = \BA(323123), \BA_8 = \BA(321323),      
$$
and the ones in the lower way by 
$$
\BB_0 = \BA_1 = \BA(121321), \BB_1 = \BA(123121), \BB_2 = \BA(123212), \BB_3 = \BA(132312), 
$$
$$
\BB_4 = \BA(132132), \BB_5 = \BA(312132), \BB_6 = \BA(321232), \BB_7 = \BA_8 = \BA(321323). 
$$

\

(a) {\em The left cell (upper road)} 

\

\subsection{Upper way.} Denote for brevity the coordinates in 
$\BA_1, \ldots, \BA_8$ by  
$$
(a_1,\ldots, c_6), (a'_1,\ldots, c'_6), (a''_1,\ldots, c''_6), (a'''_1,\ldots, c'''_6),
$$
$$ 
(a^{iv}_1,\ldots, c^{iv}_6), (a^{v}_1,\ldots, c^{v}_6), (a^{vi}_1,\ldots, c^{vi}_6), 
(a^{vii}_1,\ldots, c^{vii}_6) 
$$
respectively.

\subsection{$\bR(1)$.}\label{bR1} Let us start with
$$
\BA_1\overset{\bR(1)}\lra \BA_2.
$$
The correspondence
$$
\bR(1) = R\times\text{Id} \subset \BA_1\times \BA_2 
$$ 
is given by $6 + 9 = 15$ equations: six equations 
(\ref{RC1}) - (\ref{RC6}) on $9$ variables $a_1, \ldots, c'_3$, and $9$ equations on the remaining $9$ variables $a_4, \ldots, c'_6$:
\begin{equation}
a_4 = a'_4, b_4 = b'_4, \ldots c_6 = c'_6.
\end{equation}
So $\bR(1)$ is a cell of dimension $21 = 36 - 15$. The free coordinates in it are
$$
a_1, \ldots, c_6, a'_1, a'_3, b'_1,
$$
and the remaining variables are being
\begin{equation}
c'_1 = a'_1c_2c_3/(a_3c_1)
\end{equation}
\begin{equation}
a'_2 = a_1a_3 + a_2b_1c_3
\end{equation}
\begin{equation}
b'_2 = a_1b_3/a'_3
\end{equation}
\begin{equation}
c'_2 = a_3c_1/a'_1
\end{equation}
\begin{equation}
b'_3 = a'_3b_1b_2/(a_1b_3)
\end{equation}
\begin{equation}
c'_3 = b_3c_1/b'_1
\end{equation}
and
\begin{equation}
a'_4 = a_4, \ldots, c'_6 = c_6
\end{equation}

\subsection{$\bR(3)$.}\label{bR3} Next, 
$$
\BA_2\overset{\bR(3)}\lra \BA_3, \bR(3)\subset \BA_2\times \BA_3.
$$
It is a cell of dimension $21$.
Free coordinates in it are
$$
a'_1, \ldots, c'_6, a''_3, a''_5, b''_3,
$$
and it acts inside the third, fourth, and fifth triples. 

So the remaining coordinates are being
\begin{equation}
a''_1 = a'_1, \ldots, c''_2 = c'_2
\end{equation}
and
\begin{equation}
c''_3 = a''_3c'_4c'_5/(a'_5c'_3)
\end{equation}
\begin{equation}
a''_4 = a'_3a'_5 + a'_4b'_3c'_5
\end{equation}
\begin{equation}
b''_4 = a'_3b'_5/a''_5
\end{equation}
\begin{equation}
c''_4 = a'_5c'_3/a''_3
\end{equation}
\begin{equation}
b''_5 = a''_5b'_3b'_4/(a'_3b'_5)
\end{equation}
\begin{equation}
c''_5 = b'_5c'_3/b''_3
\end{equation}
and
\begin{equation}
a''_6 = a'_6, b''_6 = b'_6, c''_6 = c'_6
\end{equation}

\subsection{$\bR(3)\bR(1)$.}
Now let us compute the composition
$$
\BA(121321)\overset{\bR(1)}\lra \BA(212321)\overset{\bR(3)}\lra \BA(213231).
$$
We start with the fiber product 
$$
R' : = \bR(1)\times_{\BA_2}\bR(3)\subset\BA_1\times\BA_2\times \BA_3.
$$
It is given by the equations
\begin{equation}\label{R1}
A(a_1,b_1,c_1)_1 A(a_2,b_2,c_2)_2 A(a_3,b_3,c_3)_1 = 
A(a'_1,b'_1,c'_1)_2 A(a'_2,b'_2,c'_2)_1 A(a'_3,b'_3,c'_3)_2.
\end{equation}
and
\begin{equation}\label{R3}
A(a'_3,b'_3,c'_3)_2 A(a'_4,b'_4,c'_4)_3 A(a'_5,b'_5,c'_5)_2 = 
A(a''_3,b''_3,c''_3)_3 A(a''_4,b''_4,c''_4)_2 A(a''_5,b''_5,c''_5)_3
\end{equation}
which means $12$ equations:
\begin{equation}\label{RC11}
a_1a_3 + a_2b_1c_3 = a'_2
\end{equation}
\begin{equation}\label{RC12}
a_1b_3 = b'_2a'_3
\end{equation}
\begin{equation}\label{RC13}
b_1b_2 = b'_2b'_3
\end{equation}
\begin{equation}\label{RC14}
c_1a_3 = a'_1c'_2
\end{equation}
\begin{equation}\label{RC15}
c_2c_3 = c'_1c'_2
\end{equation}
\begin{equation}\label{RC16}
c_1b_3 = b'_1c'_3,
\end{equation}
and
\begin{equation}\label{RC21}
a'_3a'_5 + a'_4b'_3c'_5 = a''_4
\end{equation}
\begin{equation}\label{RC22}
a'_3b'_5 = b''_4a''_5
\end{equation}
\begin{equation}\label{RC23}
b'_3b'_4 = b''_4b''_5
\end{equation}
\begin{equation}\label{RC24}
c'_3a'_5 = a''_3c''_4
\end{equation}
\begin{equation}\label{RC25}
c'_4c'_5 = c''_3c''_4
\end{equation}
\begin{equation}\label{RC26}
c'_3b'_5 = b''_3c''_5
\end{equation}
together with $18$ simple equations
\begin{equation}
a_4 = a'_4, \ldots, c_6 = c'_6;
\end{equation}
and
\begin{equation}
a'_1 = a''_1, \ldots, c'_2 = c''_2, \ a'_6 = a''_6, b'_6 = b''_6, c'_6 = c''_6  
\end{equation}
which makes together $12 + 18 = 30$ equations, so we would expect 
$R'$ to be of dimension $54 - 30 = 24$. 

Let us see it using our parametrizations. The free variables (coordinates) in $R'$ will be 
$$
a_1, \ldots, c_6, a'_1, a'_3, b'_1, a''_3, a''_5, b''_3
$$ 
which makes $18 + 3 + 3 = 24$. 
So indeed, $R'$ is a cell of dimension $24$.  

By definition
$$
\bR(3)\bR(1) = p_{13}(\bR(1)\times_{\BA_2}\bR(3))\subset \BA_1\times\BA_3.
$$
Taking the projection $p_{13}$ we deduce that  
$\bR(3)\bR(1)$ will be a cell of dimension $24 = 18 + 3 + 3$ 
with coordinates $$
a_1, \ldots, c_6, a'_1, a'_3, b'_1, a''_3, a''_5, b''_3.
$$

This is similar to $\bR(1)$ but the formulas for the remaining variables will be 
more complicated. Namely:
\begin{equation}
a_1'' = a_1',\ 
b_1'' = b_1',\ 
c_1'' = \frac{a_1'c_2c_3}{a_3c_1} 
\end{equation}
\begin{equation}
a_2'' = a_1a_3 + a_2b_1c_3
\end{equation} 
\begin{equation}
b_2'' = \frac{a_1b_3}{a_3'},\ 
c_2'' = \frac{a_3c_1}{a_1'}
\end{equation}
\begin{equation}
c_3'' = \frac{a_3''b_1'c_4c_5}{a_5b_3c_1}
\end{equation}
\begin{equation}
a_4'' = a_3'a_5 + \frac{a_3'a_4b_1b_2c_5}{a_1b_3}
\end{equation}
\begin{equation}
b_4'' = \frac{a_3'b_5}{a_5''},\ c_4'' = \frac{a_5b_3c_1}{a_3''b_1'}
\end{equation}
\begin{equation}
b_5'' = \frac{a_5''b_1b_2b_4}{a_1b_3b_5},\  
c_5'' = \frac{b_3b_5c_1}{b_1'b_3''}, 
\end{equation}
\begin{equation}
a_6'' = a_6,\  b_6'' = b_6,\  c_6'' = c_6
\end{equation}
Summarizing,
$$
(a_1'', \ldots, c_6'') = (a_1', b_1', 
\frac{a_1'c_2c_3}{a_3c_1},  
$$
$$
a_1a_3 + a_2b_1c_3, 
\frac{a_1b_3}{a_3'}, 
\frac{a_3c_1}{a_1'},
a_3'', b_3'', \frac{a_3''b_1'c_4c_5}{a_5b_3c_1}, 
$$
$$
a_3'\biggl(a_5 + \frac{a_4b_1b_2c_5}{a_1b_3}\biggr), 
\frac{a_3'b_5}{a_5''}, \frac{a_5b_3c_1}{a_3''b_1'},
$$   
\begin{equation}\label{R3R1-resume}
a_5'', \frac{a_5''b_1b_2b_4}{a_1b_3b_5},  
\frac{b_3b_5c_1}{b_1'b_3''}, 
a_6, b_6, c_6)  
\end{equation}

\subsection{$\bL(5)\bL(2)$.} The morphism
$$
\bL(5)\bL(2): \BA_3 \lra \BA_5
$$
interchanges the second triple with the third one, and the fifth triple with the sixth one.
In other words
$$
\bL(5)\bL(2)\subset \BA_3 \times \BA_5
$$
is a closed cell of dimension $18$ with free coordinates $(a_1'', \ldots, c_6'')$, the other 
variables being 
$$
a_1^{iv} = a_1'',\ b_1^{iv} = b_1'',\ c_1^{iv} = c_1'', 
$$
$$
a_2^{iv} = a_3'',\ b_2^{iv} = b_3'',\ c_2^{iv} = c_3'',\ a_3^{iv} = a_2'',\ b_3^{iv} = b_2'',\ c_3^{iv} = c_2'', 
$$
$$
a_4^{iv} = a_4'',\ b_4^{iv} = b_4'',\ c_4^{iv} = c_4'',
$$
$$
a_5^{iv} = a_6'',\ b_5^{iv} = b_6'',\ c_5^{iv} = c_6'',\ a_6^{iv} = a_5'',\ b_6^{iv} = b_5'',\ c_6^{iv} = c_5'' 
$$ 

\subsection{$\bL(5)\bL(2)\bR(3)\bR(1)$.} 
Therefore the correspondence 
$$
\bL(5)\bL(2)\bR(3)\bR(1)\subset \BA_1\times \BA_5
$$
will be a $24$-dimensional cell with the same free coordinates as for $\bR(3)\bR(1)$, i.e. 
$$
a_1, \ldots, c_6, a'_1, a'_3, b'_1, a''_3, a''_5, b''_3,
$$
and is given by    
$$
(a_1^{iv}, \ldots, c_6^{iv}) = (a_1', b_1', 
\frac{a_1'c_2c_3}{a_3c_1},  
$$
$$
a_3'', b_3'', \frac{a_3''b_1'c_4c_5}{a_5b_3c_1}, 
a_1a_3 + a_2b_1c_3, 
\frac{a_1b_3}{a_3'}, 
\frac{a_3c_1}{a_1'}, 
$$
$$
a_3'\biggl(a_5 + \frac{a_4b_1b_2c_5}{a_1b_3}\biggr), 
\frac{a_3'b_5}{a_5''}, \frac{a_5b_3c_1}{a_3''b_1'},
$$
\begin{equation}\label{L5L2R3R1}
a_6, b_6, c_6, 
a_5'', \frac{a_5''b_1b_2b_4}{a_1b_3b_5},  
\frac{b_3b_5c_1}{b_1'b_3''})  
\end{equation}
 
We leave it to an interested reader to finish the computation of both sides of the sonnet equation.

\section{A smaller, two-parameter $R$-correspondence}\label{sect-smaller}

\subsection{•} Let us start with a $c$-upper triangular matrix
$$
A(a,b) = \left(\begin{matrix}
a & b\\
b^{-1} & 0
\end{matrix}\right)
$$
Then we have
$$
A(a_1,b_1)_1A(a_2,b_2)_2A(a_3, b_3)_1 = 
\left(\begin{matrix}
a_1a_3 + b_1a_2b_3^{-1} & a_1b_3 & b_1b_2\\
b_1^{-1}a_3 & b_1^{-1}b_3 & 0\\
b_2^{-1}b_3^{-1} & 0 & 0
\end{matrix}\right)
$$
On the other hand 
$$
A(a_1,b_1)_2A(a_2,b_2)_1A(a_3, b_3)_2 = 
\left(\begin{matrix}
a_2 & b_2a_3 & b_2b_3\\
a_1b_2^{-1} & b_1b_3^{-1} & 0\\
b_1^{-1}b_2^{-1} & 0 & 0
\end{matrix}\right)
$$
Equating
$$
A(a_1,b_1)_1A(a_2,b_2)_2A(a_3, b_3)_1 = A(a'_1,b'_1)_2A(a'_2,b'_2)_1A(a'_3, b'_3)_2
$$
we arrive at $6$ equations
$$
a_1a_3 + b_1a_2b_3^{-1} = a_2'
$$
$$
a_1b_3 = b'_2a'_3
$$
$$
b_1b_2 = b'_2b'_3
$$
$$
b_1^{-1}a_3 = a'_1b_2^{\prime -1}
$$
$$
b_1^{-1}b_3 = b'_1b_3^{\prime -1}
$$
$$
b_2^{-1}b_3^{-1} = b_1^{\prime -1}b_2^{\prime -1}
$$
In the last five equations if we fix $a_1'$, the other primed quantities are expressed 
through non-primed and $a_1'$:
$$
b_1' = \frac{a_3b_2b_3}{a_1'b_1},\
b_2' = \frac{a_1'b_1}{a_3}, \ 
a_3' = \frac{a_1a_3b_3}{a_1'b_1},\  
b_3' = \frac{a_3b_2}{a_1'}
$$
So our correspondence acts as follows:
$$
R:\ (a_1, b_1, a_2, b_2, a_3, b_3) \mapsto 
(a'_1, \frac{a_3b_2b_3}{a_1'b_1}, \frac{a_1a_3b_3 + a_2b_1}{b_3}, 
\frac{a_1'b_1}{a_3}, \frac{a_1a_3b_3}{a_1'b_1},  
\frac{a_3b_2}{a_1'}), 
$$
in other words
\begin{equation}\label{smaller-r}
(a'_1, b'_1, a'_2, b'_2, a'_3, b'_3) = 
(a'_1, \frac{a_3b_2b_3}{a_1'b_1}, \frac{a_1a_3b_3 + a_2b_1}{b_3}, 
\frac{a_1'b_1}{a_3}, \frac{a_1a_3b_3}{a_1'b_1},  
\frac{a_3b_2}{a_1'}).
\end{equation}

\subsection{•} Now we can go through the tetrahedron - sonnet diagram (\ref{tetr-eq-diagram}) and describe the corresponding $f_+, f_-$ varieties.

Each $\BA_i, \BB_j$ is isomorphic to $\BA^{12}$.
The coordinates in $\BA_i, \BB_j$ will be 
$$
(a_1, b_1, \ldots, a_6, b_6), (a'_1, b'_1 \ldots, a'_6, b'_6), 
$$
etc.

\

(a) {\em Upper way}

\

\subsection{$\bR(1)$.} For
$$
\bR(1):\ \BA_1 \lra \BA_2
$$
the coordinates will be
$$
a_1, \ldots, b_6, a_1',
$$
and it acts as follows
$$
(a'_1, b'_1, a'_2, b'_2, a'_3, b'_3) = 
(a'_1, \frac{a_3b_2b_3}{a_1'b_1}, \frac{a_1a_3b_3 + a_2b_1}{b_3}, 
\frac{a_1'b_1}{a_3}, \frac{a_1a_3b_3}{a_1'b_1},  
\frac{a_3b_2}{a_1'}),
$$
$$
(a'_4, \ldots, b'_6) = (a_4, \ldots, b_6)
$$

\subsection{$\bR(3)$.} 

$$
\bR(3):\ \BA_2 \lra \BA_3
$$
acts as follows
$$
(a''_1, \ldots, b''_2) = (a'_1, \ldots, b'_2)
$$
$$
(a''_3, b''_3, a''_4, b''_4, a''_5, b''_5) = 
(a''_3, \frac{a'_5b'_4b'_5}{a_3''b'_3}, \frac{a'_3a'_5b'_5 + a'_4b'_3}{b'_5}, 
\frac{a_3''b'_3}{a'_5}, \frac{a_3'a_5'b_5'}{a_3''b_3'},  
\frac{a'_5b'_4}{a_3''}),
$$
$$
(a_6'', b_6'') = (a_6', b_6')
$$

\subsection{$\bR(3)\bR(1)$.} Consider the composition 
$$
\bR(3)\bR(1):\ \BA_1 \lra \BA_2 \lra \BA_3
$$
The coordinates for it will be
$$
a_1, \ldots, b_6, a'_1, a''_3
$$
It acts as follows
$$
(a_1'', \ldots, b_6'') = (a_1', \frac{a_3b_2b_3}{a_1'b_1}, 
\frac{a_2b_1 + a_1a_3b_3}{b_3}, \frac{a_1'b_1}{a_3}, 
$$
$$
a_3'', \frac{a_1'a_5b_4b_5}{a_3''a_3b_2}, 
\frac{a_3(a_1a_5b_3b_5 + a_4b_1b_2)}{a_1'b_1b_5}, 
\frac{a_3''a_3b_2}{a_1'a_5},
$$
$$ 
\frac{a_1a_5b_3b_5}{a_3''b_1b_2}, \frac{a_5b_4}{a_3''}, a_6, b_6) 
$$

\subsection{$\bL(5)\bL(2)$.} The composition
$$
\bL(5)\bL(2):\ \BA_3 \lra \BA_5
$$
acts as follows
$$
(a_1^{iv}, \ldots, b_6^{iv}) = (a_1'', b_1'', a_3'', b_3'', a_2'', b_2'', 
a_4'', b_4'', a_6'', b_6'', a_5'', b_5'').
$$

\subsection{$\bL(5)\bL(2)\bR(3)\bR(1)$.} For  
$$
\bL(5)\bL(2)\bR(3)\bR(1): \BA_1 \lra \BA_5
$$
the coordinates will be 
$$
a_1, \ldots, b_6, a_1', a_3''
$$
It is acting as follows:
$$
(a_1^{iv}, \ldots, b_6^{iv}) = (a_1', \frac{a_3b_2b_3}{a_1'b_1}, 
a_3'', \frac{a_1'a_5b_4b_5}{a_3''a_3b_2}, 
$$
$$
\frac{a_2b_1 + a_1a_3b_3}{b_3}, \frac{a_1'b_1}{a_3}, 
\frac{a_3(a_1a_5b_3b_5 + a_4b_1b_2)}{a_1'b_1b_5},  
\frac{a_3''a_3b_2}{a_1'a_5},
$$
$$
a_6, b_6, \frac{a_1a_5b_3b_5}{a_3''b_1b_2}, \frac{a_5b_4}{a_3''})
$$

\subsection{$\bR(3)$} For 
$$
\bR(3):\ \BA_5 \lra \BA_6
$$
the coordinates will be 
$$
a_1^{iv}, \ldots, b_6^{iv}, a_3^v
$$
and it acts as follows:
$$
(a^v_1, \ldots, b^v_2) = (a_1^{iv}, \ldots, b_2^{iv})
$$
$$
(a_3^v, b_3^v, a_4^v, b_4^v, a_5^v, b_5^v) = 
(a_3^v, \frac{a_5^{iv}b_4^{iv}b_5^{iv}}{a_3^vb_3^{iv}}, \frac{a_3^{iv}a_5^{iv}b_5^{iv} + a_4^{iv}b_3^{iv}}{b_5^{iv}}, 
\frac{a_3^vb_3^{iv}}{a_5^{iv}}, \frac{a_3^{iv}a_5^{iv}b_5^{iv}}{a_3^vb_3^{iv}},  
\frac{a_5^{iv}b_4^{iv}}{a_3^v}),
$$
$$
(a_6^v, b_6^v) = (a_6^{iv}, b_6^{iv})
$$

\subsection{$\bR(1)$.} For 
$$
\bR(1):\ \BA_6 \lra \BA_7
$$
the coordinates will be
$$
a_1^v, \ldots, b_6^v, a_1^{vi},
$$
and it acts as follows
$$
(a_1^{vi}, b_1^{vi}, a_2^{vi}, b_2^{vi}, a_3^{vi}, b_3^{vi}) = 
(a_1^{vi}, \frac{a_3^vb_2^vb_3^v}{a_1^{vi}b_1^v}, \frac{a_1^va_3^vb_3^v + a_2^vb_1^v}{b_3^v}, 
\frac{a_1^{vi}b_1^v}{a_3^v}, \frac{a_1^va_3^vb_3^v}{a_1^{vi}b_1^v},  
\frac{a_3^vb_2^v}{a_1^{vi}}),
$$
$$
(a_4^{vi}, \ldots, b_6^{vi}) = (a_4^v, \ldots, b_6^v)
$$ 

\subsection{$\bR(1)\bR(3)$.} For 
$$
\bR(1)\bR(3): \BA_5 \lra \BA_7
$$
the coordinates will be
$$
a_1^{iv}, \ldots, b_6^{iv}, a_3^v, a_1^{vi}, 
$$
and it acts as follows
$$
(a_1^{vi}, \ldots, b_6^{vi}) = (a_1^{vi}, 
\frac{a_5^{iv}b_2^{iv}b_4^{iv}b_5^{iv}}{a_1^{vi}b_1^{iv}b_3^{iv}}, 
\frac{a_3^v(a_2^{iv}b_1^{iv}b_3^{iv} + a_1^{iv}a_5^{iv}b_4^{iv}b_5^{iv})}{a_5^{iv}b_4^{iv}b_6^{iv}}, 
\frac{a_1^{vi}b_1^{iv}}{a_3^v}, 
$$
$$
\frac{a_1^{iv}a_5^{iv}b_4^{iv}b_5^{iv}}{a_1^{vi}b_1^{iv}b_3^{iv}},
\frac{a_3^vb_2^{iv}}{a_1^{vi}}, 
\frac{a_4^{iv}b_3^{iv} + a_3^{iv}a_5^{iv}b_5^{iv}}{b_5^{iv}}, 
\frac{a_3^{v}b_3^{iv}}{a_5^{iv}}, 
$$
$$
\frac{a_3^{iv}a_5^{iv}b_5^{iv}}{a_3^{v}b_3^{iv}},
\frac{a_5^{iv}b_4^{iv}}{a_3^{v}}, a_6^{iv}, b_6^{iv})  
$$

\subsection{$\bL(3)$.} The map
$$
\bL(3):\ \BA_7 \lra \BA_8
$$
acts as follows
$$
(a_1^{vii}, b_1^{vii}, a_2^{vii}, b_2^{vii}, a_3^{vii}, b_3^{vii},  
a_4^{vii}, b_4^{vii}, a_5^{vii}, b_5^{vii}, a_6^{vii}, b_6^{vii}) = 
$$
$$
(a_1^{vi}, b_1^{vi}, a_2^{vi}, b_2^{vi}, a_4^{vi}, b_4^{vi},  
a_3^{vi}, b_3^{vi}, a_5^{vi}, b_5^{vi}, a_6^{vi}, b_6^{vi})  
$$
(we just interchanged the third group with the fourth one).

\subsection{$\bL(3)\bR(1)\bR(3)$.} For the map 
$$
\bL(3)\bR(1)\bR(3): \BA_5 \lra \BA_8
$$
the coordinates will be
$$
a_1^{iv}, \ldots, b_6^{iv}, a_3^v, a_1^{vi}, 
$$
and it acts as follows
$$
(a_1^{vii}, \ldots, b_6^{vii}) = (a_1^{vi}, 
\frac{a_5^{iv}b_2^{iv}b_4^{iv}b_5^{iv}}{a_1^{vi}b_1^{iv}b_3^{iv}}, 
\frac{a_3^v(a_2^{iv}b_1^{iv}b_3^{iv} + a_1^{iv}a_5^{iv}b_4^{iv}b_5^{iv})}{a_5^{iv}b_4^{iv}b_6^{iv}}, 
\frac{a_1^{vi}b_1^{iv}}{a_3^v}, 
$$
$$
\frac{a_4^{iv}b_3^{iv} + a_3^{iv}a_5^{iv}b_5^{iv}}{b_5^{iv}}, 
\frac{a_3^{v}b_3^{iv}}{a_5^{iv}}, 
\frac{a_1^{iv}a_5^{iv}b_4^{iv}b_5^{iv}}{a_1^{vi}b_1^{iv}b_3^{iv}},
\frac{a_3^vb_2^{iv}}{a_1^{vi}}, 
$$
$$
\frac{a_3^{iv}a_5^{iv}b_5^{iv}}{a_3^{v}b_3^{iv}},
\frac{a_5^{iv}b_4^{iv}}{a_3^{v}}, a_6^{iv}, b_6^{iv})  
$$

\subsection{LHS.} The left hand side of the tetrahedron equation will be
$$
C_+^m:= \bL(3)\bR(1)\bR(3)\bL(5)\bL(2)\bR(3)\bR(1):\ \BA_1 \lra \BA_8
$$
The coordinates are
$$
a_1, \ldots, b_6, a_1', a_3'', a_3^v, a_1^{vi}, 
$$
and it acts as follows
$$
(a_1^{vii}, \ldots, b_6^{vii}) = (a_1^{vi}, 
\frac{a_6b_4b_5b_6}{a_1^{vi}b_2b_3}, 
\frac{a_1'a_3''a_3^v(a_5b_3 + a_3a_6b_6)}{a_3a_5a_6b_4},  
\frac{a_1^{vi}a_3b_2b_3}{a_1'a_3^vb_1},  
$$
$$
\frac{a_1a_5b_3^2b_5 + a_4b_1b_2b_3 + (a_2b_1 + a_1a_3b_3)a_6b_6b_5}{b_3b_6b_5},
\frac{a_1'a_3^vb_1}{a_3a_6}, 
\frac{a_3''a_3a_6b_6}{a_1^{vi}a_5b_3}, 
\frac{a_1'a_3^va_5b_4b_5}{a_3''a_1^{vi}a_3b_2}, 
$$
$$
\frac{(a_2b_1 + a_1a_3b_3)a_3a_6b_6}{a_1'a_3^vb_1b_3},
\frac{a_3''a_3a_6b_2}{a_1'a_3^va_5},
\frac{a_1a_5b_3b_5}{a_3''b_1b_2}, 
\frac{a_5b_4}{a_3''})
$$

\

(b) {\em Lower way}

\

\subsection{$\bL(3)$.} The map
$$
\bL(3):\ \BB_0 \lra \BB_1
$$
acts as follows
$$
(a'_1, \ldots, b'_2) = (a_1, \ldots, b_2),
$$
$$
(a'_3, b'_3, a'_4, b'_4) = (a_4, b_4, a_3, b_3), 
$$
$$
(a'_5, \ldots, b'_6) = (a_5, \ldots, b_6)
$$

\subsection{$\bR(4)$.} The map
$$
\bR(4):\ \BB_1 \lra \BB_2
$$
acts as follows
$$
(a_1'', \ldots, b_3'') = (a_1', \ldots, b_3'), 
$$
$$
(a_4'', \ldots, b_6'') = (a_4'', \frac{a_6'b_5'b_6'}{a_4''b_4'}, 
\frac{a_5'b_4' + a_4'a_6'b_6'}{b_6'}, \frac{a_4''b_4'}{a_6'}, 
\frac{a_4'a_6'b_6'}{a_4''b_4'}, \frac{a_6'b_5'}{a_4''})
$$

\subsection{$\bR(2)$.} The map
$$
\bR(2):\ \BB_2 \lra \BB_3
$$
acts as follows:
$$
(a_1''', b_1''') = (a_1'', b_1''),
$$
$$
(a_2''', \ldots, b_4''') = (a_2''', \frac{a_4''b_3''b_4''}{a_2'''b_2''}, 
\frac{a_3''b_2'' + a_2''a_4''b_4''}{b_4''}, \frac{a_2'''b_2''}{a_4''}, 
\frac{a_2''a_4''b_4''}{a_2'''b_2''}, \frac{a_4''b_3''}{a_2'''}), 
$$
$$
(a_5''', \ldots, b_6''') = (a_5'', \ldots, b_6'')
$$

\subsection{$\bR(2)\bR(4)$.} The map
$$
\bR(2)\bR(4):\ \BB_1 \lra \BB_3
$$
has the coordinates
$$
(a_1', \ldots, b_6', a_2''', a_4''),
$$
and acts as follows:
$$
(a_1''', \ldots, b_6''') = (a_1', b_1', a_2''', \frac{a_6'b_3'b_5'b_6'}{a_2'''b_2'b_4'}, 
$$
$$
\frac{a_4''(a_3'b_2'b_4' + a_2'a_6'b_5'b_6')}{a_6'b_5'b_6'},
\frac{a_2'''b_2'}{a_4''}, 
\frac{a_2'a_6'b_5'b_6'}{a_2'''b_2'b_4'}, 
\frac{a_4''b_3'}{a_2'''}, 
$$
$$
\frac{a_5'b_4' + a_4'a_6'b_6'}{b_6'},
\frac{a_4''b_4'}{a_6'}, 
\frac{a_4'a_6'b_6'}{a_4''b_4'}, 
\frac{a_6'b_5'}{a_4''}) 
$$

\subsection{$C^m_{-+} = \bR(2)\bR(4)\bL(3)$.} The map
$$
C^m_{-+} := \bR(2)\bR(4)\bL(3):\ \BB_0 \lra \BB_3
$$
has the coordinates
$$
(a_1, \ldots, b_6, a_2''', a_4''),
$$
and acts as follows:
$$
(a_1''', \ldots, b_6''') = (a_1, b_1, a_2''', \frac{a_6b_4b_5b_6}{a_2'''b_2b_3}, 
$$
$$
\frac{a_4''(a_4b_2b_3 + a_2a_6b_5b_6)}{a_6b_5b_6},
\frac{a_2'''b_2}{a_4''}, 
\frac{a_2a_6b_5b_6}{a_2'''b_2b_3}, 
\frac{a_4''b_4}{a_2'''}, 
$$
$$
\frac{a_5b_3 + a_3a_6b_6}{b_6},
\frac{a_4''b_3}{a_6}, 
\frac{a_3a_6b_6}{a_4''b_3}, 
\frac{a_6b_5}{a_4''}) 
$$

\


\subsection{$\bL(1)\bL(4)$.} 
The map 
$$
\bL(1)\bL(4): \BB_3 \lra \BB_5
$$
acts as follows:
$$
(a_1^v, b_1^v, a_2^v, b_2^v) = (a_2''', b_2''', a_1''', b_1''')
$$
$$
(a_3^v, b_3^v) = (a_3''', b_3''')
$$
$$
(a_4^v, b_4^v, a_5^v, b_5^v) = (a_5''', b_5''', a_4''', b_4''')
$$
$$
(a_6^v, b_6^v) = (a_6''', b_6''')
$$

\subsection{$\bR(2)$.} The map 
$$
\bR(2): \BB_5 \lra \BB_6
$$
acts as follows:
$$
(a_1^{vi}, b_1^{vi}) = (a_1^{v}, b_1^{v})
$$
$$
(a_2^{vi}, \ldots, b_4^{vi}) = (a_2^{vi}, \frac{a_4^{v}b_3^{v}b_4^{v}}{a_2^{vi}b_2^{v}}, 
\frac{a_3^{v}b_2^{v} + a_2^{v}a_4^{v}b_4^{v}}{b_4^{v}}, \frac{a_2^{vi}b_2^{v}}{a_4^{v}}, 
\frac{a_2^{v}a_4^{v}b_4^{v}}{a_2^{vi}b_2^{v}}, \frac{a_4^{v}b_3^{v}}{a_2^{vi}}) 
$$
$$
(a_5^{vi}, \ldots, b_6^{vi}) = (a_5^{v}, \ldots, b_6^{v})
$$

\subsection{$\bR(4)$.} The map 
$$
\bR(4): \BB_6 \lra \BB_7
$$
acts as follows:
$$
(a_1^{vii}, \ldots, b_3^{vii}) = (a_1^{vi}, \ldots, b_3^{vi})
$$
$$
(a_4^{vii}, \ldots, b_6^{vii}) = (a_4^{vii}, \frac{a_6^{vi}b_5^{vi}b_6^{vi}}{a_4^{vii}b_4^{vi}}, 
\frac{a_5^{vi}b_4^{vi} + a_4^{vi}a_6^{vi}b_6^{vi}}{b_6^{vi}}, \frac{a_4^{vii}b_4^{vi}}{a_6^{vi}}, 
\frac{a_4^{vi}a_6^{vi}b_6^{vi}}{a_4^{vii}b_4^{vi}}, \frac{a_6^{vi}b_5^{vi}}{a_4^{vii}}) 
$$

\subsection{$\bR(4)\bR(2)$.} The map 
$$
\bR(4)\bR(2): \BB_5 \lra \BB_7
$$
acts as follows:
$$
(a_1^{vii}, \ldots, b_6^{vii}) = (a_1^v, b_1^v, a_2^{vi}, \frac{a_4^vb_3^vb_4^v}{a_2^{vi}b_2^v}, 
$$
$$
\frac{a_3^vb_2^v + a_2^va_4^vb_4^v}{b_4^v}, 
\frac{a_2^{vi}b_2^v}{a_4^v}, 
a_4^{vii},
\frac{a_2^va_6^vb_5^vb_6^v}{a_4^{vii}a_4^vb_3^v}, 
$$
$$
\frac{a_4^v(a_5^vb_2^vb_3^v + a_2^va_6^vb_4^vb_6^v)}{a_2^{vi}b_2^vb_6^v},
\frac{a_4^{vii}a_4^vb_3^v}{a_2^{vi}a_6^v}, 
\frac{a_2^va_6^vb_4^vb_6^v}{a_4^{vii}b_2^vb_3^v},
\frac{a_6^vb_5^v}{a_4^{vii}})
$$

\subsection{$C^m_{--} = \bR(4)\bR(2)\bL(1)\bL(4)$.} The map 
$$
C^m_{--} : =\bR(4)\bR(2)\bL(1)\bL(4): \BB_3 \lra \BB_7
$$
is acting as follows:
$$
(a_1^{vii}, \ldots, b_6^{vii}) = (a_2''', b_2''', a_2^{vi}, 
\frac{a_5'''b_3'''b_5'''}{a_2^{vi}b_1'''}, 
$$
$$
\frac{a_3'''b_1''' + a_1'''a_5'''b_5'''}{b_5'''}, 
\frac{a_2^{vi}b_1'''}{a_5'''}, 
a_4^{vii},
\frac{a_1'''a_6'''b_4'''b_6'''}{a_4^{vii}a_5'''b_3'''}, 
$$
$$
\frac{a_5'''(a_4'''b_1'''b_3''' + a_1'''a_6'''b_5'''b_6''')}{a_2^{vi}b_1'''b_6'''},
\frac{a_4^{vii}a_5'''b_3'''}{a_2^{vi}a_6'''}, 
\frac{a_1'''a_6'''b_5'''b_6'''}{a_4^{vii}b_1'''b_3'''},
\frac{a_6'''b_4'''}{a_4^{vii}})
$$

\subsection{RHS.} The right hand side of the tetrahedron equation will be
$$
C_-^m:= C^m_{--}C^m_{-+} = \bR(4)\bR(2)\bL(1)\bL(4)\bR(2)\bR(4)\bL(3) :\ \BB_0 \lra \BB_7
$$
The coordinates are
$$
a_1, \ldots, b_6, a_4'', a_2''', a_2^{vi}, a_4^{vii}, 
$$
and it acts as follows
$$
(a_1^{vii}, \ldots, b_6^{vii}) = (a_2''', b_1, a_2^{vi}, 
\frac{a_2'''(a_5b_3 + a_3a_6b_6)b_2b_3}{a_2^{vi}a_6b_1b_6}, 
$$
$$
\frac{a_4b_1b_2b_3 + a_1a_5b_3^2b_5 + a_2a_6b_1b_5b_6 + a_1a_3a_6b_3b_5b_6}{b_3b_5b_6}, 
\frac{a_2^{vi}b_1b_6}{a_5b_3 + a_3a_6b_6}, 
a_4^{vii},
\frac{a_1a_3a_6^2b_4b_5b_6^2}{a_2'''a_4^{vii 2}b_2b_3(a_5b_3 + a_3a_6b_6)}, 
$$
$$
\frac{(a_5b_3 + a_3a_6b_6)(a_2b_1 + a_1a_3b_3)}{a_2^{vi}b_1b_3}, 
$$
$$
\frac{a_2'''a_4^{vii}(a_5b_3 + a_3a_6b_6)b_2b_3}{a_2^{vi}a_3a_6b_6^2}, 
\frac{a_1a_3a_6b_5b_6}{a_2'''a_4^{vii}b_1b_2},
\frac{a_3a_6b_4b_6}{a_2'''a_4^{vii}b_3})
$$

\subsection{Tetrahedron equations.} The coordinates:
$$
a_1, \ldots, b_6; a_1', a_3'', a_3^v, a_1^{vi}
$$
for $C_+^m$, and
$$
a_1, \ldots, b_6; a_4'', a_2''', a_2^{vi}, a_4^{vii}
$$
for $C_-^m$.  

The equations are
$$
(a_1^{vi}, 
\frac{a_6b_4b_5b_6}{a_1^{vi}b_2b_3}, 
\frac{a_1'a_3''a_3^v(a_5b_3 + a_3a_6b_6)}{a_3a_5a_6b_4},  
\frac{a_1^{vi}a_3b_2b_3}{a_1'a_3^vb_1},  
$$
$$
\frac{a_1a_5b_3^2b_5 + a_4b_1b_2b_3 + (a_2b_1 + a_1a_3b_3)a_6b_6b_5}{b_3b_6b_5},
\frac{a_1'a_3^vb_1}{a_3a_6}, 
\frac{a_3''a_3a_6b_6}{a_1^{vi}a_5b_3}, 
\frac{a_1'a_3^va_5b_4b_5}{a_3''a_1^{vi}a_3b_2}, 
$$
$$
\frac{(a_2b_1 + a_1a_3b_3)a_3a_6b_6}{a_1'a_3^vb_1b_3},
\frac{a_3''a_3a_6b_2}{a_1'a_3^va_5},
\frac{a_1a_5b_3b_5}{a_3''b_1b_2}, 
\frac{a_5b_4}{a_3''}) = 
$$ 
$$
= (a_2''', b_1, a_2^{vi}, 
\frac{a_2'''(a_5b_3 + a_3a_6b_6)b_2b_3}{a_2^{vi}a_6b_1b_6}, 
$$
$$
\frac{a_4b_1b_2b_3 + a_1a_5b_3^2b_5 + a_2a_6b_1b_5b_6 + a_1a_3a_6b_3b_5b_6}{b_3b_5b_6}, 
\frac{a_2^{vi}b_1b_6}{a_5b_3 + a_3a_6b_6}, 
a_4^{vii},
\frac{a_1a_3a_6^2b_4b_5b_6^2}{(a_2''')^2a_4^{vii}b_2b_3(a_5b_3 + a_3a_6b_6)}, 
$$
$$
\frac{(a_5b_3 + a_3a_6b_6)(a_2b_1 + a_1a_3b_3)}{a_2^{vi}b_1b_3}, 
\frac{a_2'''a_4^{vii}(a_5b_3 + a_3a_6b_6)b_2b_3}{a_2^{vi}a_3a_6b_6^2},
$$
$$ 
\frac{a_1a_3a_6b_5b_6}{a_2'''a_4^{vii}b_1b_2},
\frac{a_3a_6b_4b_6}{a_2'''a_4^{vii}b_3}) 
$$

\subsection{•} 
In other words:

$(1a)$
$$ 
a_1^{vi} = a_2'''
$$

$(1b)$
$$
b_1 = \frac{a_6b_4b_5b_6}{a_1^{vi}b_2b_3} 
$$
whence

$(1b1)$
$$
a_1^{vi} = \frac{a_6b_4b_5b_6}{b_1b_2b_3}
$$

$(2a)$
$$
\frac{a_1'a_3''a_3^v(a_5b_3 + a_3a_6b_6)}{a_3a_5a_6b_4} = a_2^{vi}
$$

$(2b)$
$$
\frac{a_1^{vi}a_3b_2b_3}{a_1'a_3^vb_1} = 
\frac{a_2'''(a_5b_3 + a_3a_6b_6)b_2b_3}{a_2^{vi}a_6b_1b_6}
$$

whence

$(2b1)$
$$
a_2^{vi} = \frac{(a_5b_3 + a_3a_6b_6)a_1'a_3^v}{a_3a_6b_6}
$$

It follows:
 
$(2b2)$
$$
a_3'' = \frac{a_5b_4}{b_6}
$$

$(3a)$ and $(3b)$ are tautological.

$(4a)$
$$
\frac{a_3''a_3a_6b_6}{a_1^{vi}a_5b_3} = a_4^{vii}
$$
whence

$(4a1)$
$$
a_4^{vii} = \frac{a_3b_1b_2}{b_5b_6}
$$

$(4b)$
$$
\frac{a_1'a_3^va_5b_4b_5}{a_3''a_1^{vi}a_3b_2} = 
\frac{a_1a_3a_6^2b_4b_5b_6^2}{(a_2''')^2a_4^{vii}b_2b_3(a_5b_3 + a_3a_6b_6)}
$$
whence

$(4b1)$
$$
a_1'a_3^v = \frac{a_1a_3a_6b^2_6}{b_4(a_5b_3 + a_3a_6b_6)}.
$$
It follows from $(2b1)$ that

$(4b2)$
$$
a_2^{vi} = \frac{a_1b_6}{b_4}
$$

$(5a)$
$$
\frac{(a_2b_1 + a_1a_3b_3)a_3a_6b_6}{a_1'a_3^vb_1b_3} = 
\frac{(a_5b_3 + a_3a_6b_6)(a_2b_1 + a_1a_3b_3)}{a_2^{vi}b_1b_3}
$$
This coincides with $(2b1)$.

$(5b)$
$$
\frac{a_3''a_3a_6b_2}{a_1'a_3^va_5} = 
\frac{a_2'''a_4^{vii}(a_5b_3 + a_3a_6b_6)b_2b_3}{a_2^{vi}a_3a_6b_6^2}
$$
whence

$(5b1)$
$$
a_2''' = \frac{a_6b_4b_5b_6}{b_1b_2b_3},
$$
this is tautological.

$(6a)$ is tautological due to $(5b1)$

$(6b)$
$$
\frac{a_5b_4}{a_3''} = 
\frac{a_3a_6b_4b_6}{a_2'''a_4^{vii}b_3} 
$$
whence

\

$(6b1)$
$$
a_2'''a_4^{vii} = \frac{a_3a_6b_4}{b_3},
$$
therefore
$$
a_2''' = \frac{a_3a_6b_4}{a_4^{vii}b_3} = \frac{a_6b_4b_5b_6}{b_1b_2b_3},
$$
we already know that.

\

There is no condition on $a_4''$, whereas $a_1', a_3^v$ appear only through $(4b1)$.

\subsection{The sonnet equation holds true.} We conclude that the sonnet equations are satisfied, in the following sense:  
two cells 
$$
f_+ := Im(C^m_+), f_- := Im(C^m_-) \subset \BA_1 \times \BA_8 = \BB_0\times \BB_7
$$
coincide birationally, i.e. $f_+\cap f_-$ is Zarisky dense in both of them. 

Both $f_+$ and $f_-$ are birationally isomorphic to $\BA^{12}$ with coordinates 
$(a_1, b_1, \ldots, a_6, b_6$).

\

When $b\lra 1$ we get the equation discussed in the next Section  
\ref{section-very-small}.

\subsection{A quasiinverse $R$-correspondence.} Let $\BA, \BA'$ denote six dimensional 
affine spaces with coordinates $(a_1, b_1, a_2, b_2, a_3, b_3)$ and 
$(a'_1, b'_1, a'_2, b'_2, a'_3, b'_3)$ respectively. 
Define a correspondence $S:\ \BA \lra \BA'$ by     
\begin{equation}\label{inverse-to-smaller-r}
(a'_1, b'_1, a'_2, b'_2, a'_3, b'_3) = 
(a'_1, \frac{a_3b_2b_3}{a_1'b_1}, \frac{- a_1a_3b_3 + a_2b_1}{b_3}, 
\frac{a_1'b_1}{a_3}, \frac{a_1a_3b_3}{a_1'b_1},  
\frac{a_3b_2}{a_1'}).
\end{equation}
So we have changed one sign in (\ref{smaller-r}). 

\

Let $\BA''$ denote a denote six dimensional 
affine space with coordinates $(a''_1, b''_1, a''_2, b''_2, a''_3, b''_3)$. 
It is not difficult to see that the composition 
$$
SR:\ \BA\lra \BA''
$$
is given by
\begin{equation}
(a_1'', \ldots, b_3'') = (a_1'', \frac{a_1b_1}{a_1''}, \frac{a_2b_1}{b_3}, 
\frac{a_1''b_2}{a_1}, \frac{a_1a_3}{a_1''}, \frac{a_1b_3}{a_1''}).
\end{equation}
If we put here $a_1'' = a_1$ and $b_1 = b_3$ we get the identity map.

\

Moreover, if we put $b = 1$ in (\ref{inverse-to-smaller-r}) we get the inverse 
to the "very small" $R$-matrix discussed in the following Section, see 
(\ref{inv-very-small}).

\section{A very small $R$-matrix}\label{section-very-small}

Here we describe the specialization of the contents of the previous Section at $b = 1$. 
In this case $R$-correspondece turns into an $R$-matrix. 

\subsection{•} We start from a matrix 
$$
A(a) = \left(\begin{matrix}
a & 1\\
1 & 0
\end{matrix}\right)
$$
Then we have
$$
A(a_1)_1A(a_2)_2A(a_3)_1 = 
\left(\begin{matrix}
a_1a_3 + a_2 & a_1 & 1\\
a_3 & 1 & 0\\
1 & 0 & 0
\end{matrix}\right)
$$
and 
$$
A(a_1)_2A(a_2)_1A(a_3)_2 = 
\left(\begin{matrix}
a_2 & a_3 & 1\\
a_1 & 1 & 0\\
1 & 0 & 0
\end{matrix}\right),
$$
whence the $R$-matrix (not an $R$-correspondence this time)
$$
\bR:\ \BA(121) \lra \BA(212)
$$
given by
\begin{equation}\label{very-small-r}
(a_1', a_2', a_3') = (a_3, a_1a_3 + a_2, a_1)
\end{equation}
Its inverse is given by
\begin{equation}\label{inv-very-small}
(a_1', a_2', a_3') = (a_3, - a_1a_3 + a_2, a_1)
\end{equation} 

\subsection{Sonnet equations.} Let us check the tetrahedron (sonnet) equations.
Below $\BA_i = \BA^6$.

The map
$$
\bR(1): \ \BA_1 \lra \BA_2
$$
acts as follows
$$
(a'_1, \ldots , a'_6) = (a_3, a_1a_3 + a_2, a_1, a_4, a_5, a_6).
$$ 
The map
$$
\bR(3): \ \BA_2 \lra \BA_3
$$
acts as follows
$$
(a''_1, \ldots , a''_6) = (a_1', a_2', a_5', a_3'a_5' + a_4', a_3', a_6')
$$
Therefore
$$
\bR(3)\bR(1): \BA_1\lra \BA_3
$$
acts as 
$$
(a_1'', \ldots, a_6'') = (a_3, a_1a_3 + a_2, a_5, a_1a_5 + a_4, a_1, a_6)
$$

\subsection{$\bL(5)\bL(2)\bR(3)\bR(1)$} Next, 
$$
\bL(5)\bL(2): \BA_3 \lra \BA_5
$$
interchenges $2$ with $3$ and $5$ with $6$. 

It follows that 
$$
\bL(5)\bL(2)\bR(3)\bR(1): \BA_1 \lra \BA_5
$$
acts as 
$$
(a_1^{iv}, \ldots, a_6^{iv}) = (a_3, a_5, a_1a_3 + a_2, a_1a_5 + a_4, a_6, a_1)
$$

\subsection{•} The map
$$
\bR(3): \BA_5 \lra \BA_6
$$
acts as
$$
(a_1^v, \ldots, a_6^v) = (a_1^{iv}, a_2^{iv}, a_5^{iv}, a_3^{iv}a_5^{iv} + a_4^{iv}, a_3^{iv}, a_6^{iv})
$$
The map
$$
\bR(1): \ \BA_6 \lra \BA_7
$$
acts as follows
$$
(a^{vi}_1, \ldots , a^{vi}_6) = (a_3^{v}, a_1^{v}a_3^{v} + a_2^{v}, a_1^{v}, a_4^{v},  a_5^{v}, a_6^{v}).
$$
Therefore
$$
\bR(1)\bR(3): \BA_5\lra \BA_7
$$
acts as 
$$
(a_1^{vi}, \ldots, a_6^{vi}) = (a_5^{iv}, a_1^{iv}a_5^{iv} + a_2^{iv}, a_1^{iv}, 
a_3^{iv}a_5^{iv} + a_4^{iv}, a_3^{iv}, a_6^{iv})
$$

\subsection{$\bL(3)\bR(1)\bR(3)$} Next, 
$$
\bL(3): \BA_7 \lra \BA_8
$$
interchanges $3$ with $4$.

It follows that 
$$
\bL(3)\bR(1)\bR(3): \BA_5 \lra \BA_8
$$
acts as
$$
(a_1^{vii}, \ldots, a_6^{vii}) = (a_5^{iv}, a_1^{iv}a_5^{iv} + a_2^{iv}, 
a_3^{iv}a_5^{iv} + a_4^{iv}, a_1^{iv},  a_3^{iv}, a_6^{iv})
$$

\subsection{The left hand side.} The LHS of the sonnet equation is
$$
C_+^{mm} = \bL(3)\bR(1)\bR(3)\bL(5)\bL(2)\bR(3)\bR(1): \BA_1 \lra \BA_5 \lra \BA_8
$$
and acts as follows
$$
(a_1^{vii}, \ldots, a_6^{vii}) = (a_6, a_3a_6 + a_5, 
a_1a_3a_6 + a_2a_6 + a_1a_5 + a_4, a_3, 
a_1a_3 + a_2, a_1)
$$

\

{\it Lower road}

\

\subsection{•} The map 
$$
\bR(4):\ \BB_1 \lra \BB_2
$$
acts as follows
$$
(a_1'', \ldots, a_6'') = (a_1', a_2', a_3', a_6', a_4'a_6' + a_5', a_4') 
$$
The map 
$$
\bR(2):\ \BB_2 \lra \BB_3
$$
acts as follows
$$
(a_1''', \ldots, a_6''') = (a_1'', a_4'', a_2''a_4'' + a_3'', a_2'', a_5'', a_6'')
$$
Therefore
$$
\bR(2)\bR(4):\ \BB_1 \lra \BB_3
$$
acts as follows
$$
(a_1''', \ldots, a_6''') = (a_1', a_6', a_2'a_6' + a_3', 
a_2', a_4'a_6' + a_5', a_4')
$$

\subsection{•} The map 
$$
\bL(3): \BB_0 \lra \BB_1
$$
interchanges $3$ with $4$. 

Therefore the map
$$
\bR(2)\bR(4)\bL(3):\ \BB_0 \lra \BB_3
$$
acts as follows
$$
(a_1''', \ldots, a_6''') = (a_1, a_6, a_2a_6 + a_4, 
a_2, a_3a_6 + a_5, a_3)
$$

\subsection{•} The map 
$$
\bR(2): \BB_5 \lra \BB_6
$$
acts as
$$
(a_1^{vi}, \ldots, a_6^{vi}) = (a_1^v, a_4^v, a_2^va_4^v + a_3^v, a_2^v, a_5^v, a_6^v)
$$
The map 
$$
\bR(4): \BB_6 \lra \BB_7
$$
acts as
$$
(a_1^{vii}, \ldots, a_6^{vii}) = (a_1^{vi}, a_2^{vi}, a_3^{vi}, a_6^{vi}, 
a_4^{vi}a_6^{vi} + a_5^{vi}, a_4^{vi}) 
$$
Therefore
$$
\bR(4)\bR(2): \BB_5 \lra \BB_7
$$
acts as
$$
(a_1^{vii}, \ldots, a_6^{vii}) = (a_1^v, a_4^v, a_2^va_4^v + a_3^v, 
a_6^v, a_2^va_6^v + a_5^v, a_2^v)
$$

\subsection{•} The map 
$$
\bL(1)\bL(4): \BB_3 \lra \BB_5
$$
interchanges $1$ with $2$ and $4$ with $5$. 

It follows that 
$$
\bR(4)\bR(2)\bL(1)\bL(4): \BB_3 \lra \BB_7
$$
acts as
$$
(a_1^{vii}, \ldots, a_6^{vii}) = 
(a_2''', a_5''', a_1'''a_5''' + a_3''', 
a_6''', a_1'''a_6''' + a_4''', a_1''')
$$

\subsection{The right hand side.} The RHS of the sonnet equation is
$$
C_6^{mm} = \bR(4)\bR(2)\bL(1)\bL(4)\bR(2)\bR(4)\bL(3): \BB_0 \lra \BB_3 \lra \BB_7
$$
and acts as follows
$$
(a_1^{vii}, \ldots, a_6^{vii}) = 
(a_6, a_3a_6 + a_5, a_1a_3a_6 + a_1a_5 + a_2a_6 + a_4, 
a_3, a_1a_3 + a_2, a_1),
$$
which coincides with the LHS, i.e. the sonnet equations are fulfilled. 

\section{Another triple operation}

\subsection{•} For a $2\times 2$ matrix
$$
M = \left(\begin{matrix}
a & b\\
c & d
\end{matrix}\right)
$$
let us denote by $M_{13}$ the matrix 
$$
M_{13} = \left(\begin{matrix}
a & 0 & b\\
0 & 1 & 0\\
c & 0 & d
\end{matrix}\right).
$$

\subsection{Triangularity lemma.} {\em If $A, B, C$ are $c$-upper triangular 
$2\times 2$ matrices then both triple products 
$A_{12}B_{13}C_{23}$ and $A_{23}B_{13}C_{12}$ are $c$-upper triangular.} $\Box$

\subsection{$R$-matrix.} Consider a matrix  
$$
A(a) = \left(\begin{matrix}
a & 1\\
1 & 0
\end{matrix}\right).
$$
Then
$$
A(a_1)_{12}A(a_2)_{13}A(a_3)_{23} = \left(
\begin{matrix}
a_1a_2 & a_1 + a_3 & 1\\
a_2 & 1 & 0\\
1 & 0 & 0
\end{matrix}\right),
$$
and
$$
A(a_1)_{23}A(a_2)_{13}A(a_3)_{12} = \left(
\begin{matrix}
a_2a_3 & a_2 & 1\\
a_1 + a_3 & 1 & 0\\
1 & 0 & 0
\end{matrix}\right)
$$
We define the corresponding $R$-matrix by the equation
$$
A(a_1)_{12}A(a_2)_{13}A(a_3)_{23} = A(a_1')_{23}A(a_2')_{13}A(a_3')_{12}
$$
which is equivalent to a system
$$
a_1a_2 = a_2'a_3'
$$
$$
a_1 + a_3 = a_2'
$$
$$
a_2 = a_1' + a_3'
$$
or
$$
R: (a_1', a_2', a_3') = (\frac{a_2a_3}{a_1 + a_3}, a_1 + a_3, \frac{a_1a_2}{a_1 + a_3}).
$$

\subsection{$R$ is idempotent.} Note that $R^2 = \text{Id}$. 

\subsection{Sonnet equation.}
This $R$-matrix should satisfy a sonnet equation of the form
$$
L(3)R(123)R(124)L(1)L(4)R(134)R(234) = 
$$
$$
= R(234)R(134)L(2)L(5)R(124)R(123)L(3),
$$
or symbolically
$$
LRRLLRR = RRLLRRL
$$

We will use the notation:
$$
\BA(12,13,14,23,24,34) = \BA_1,\ \BA(12,13,14,34,24,23) = \BA_2
$$
$$
\BA(12,34,14,13,24,23) = \BA_3,\ \BA(34,12,14,24,13,23) = \BA_5,
$$
$$
\BA(34,24,14,12,13,23) = \BA_6,\ \BA(34,24,14,23,13,12) = \BA_7,  
$$
$$
\BA(34,24,23,14,13,12) = \BA_8
$$
for the left hand side, and
$$
\BA_1 = \BB_0,\ \BA(12,13,23,14,24,34) = \BB_1, 
$$
$$
\BA(23,13,12,14,24,34) = \BB_2,\ \BA(23,13,24,14,12,34) = \BB_3, 
$$
$$
\BA(23,24,13,14,34,12) = \BB_5, \BA(23,24,34,14,13,12) = \BB_6, 
\BA_8 = \BB_7 
$$
for the right hand side. 

\subsection{$C_{++}$.} The map
$$
R(234): \BA_1 \lra \BA_2
$$
acts as
$$
(a'_1, \ldots, a'_6) = (a_1, a_2, a_3, \frac{a_5a_6}{a_4 + a_6}, 
a_4 + a_6, \frac{a_4a_5}{a_4 + a_6})
$$
The map
$$
R(134): \BA_2 \lra \BA_3
$$
acts as
$$
(a''_1, \ldots, a''_6) = (a'_1, \frac{a'_3a'_4}{a'_2 + a'_4}, 
a'_2 + a'_4, \frac{a'_2a'_3}{a'_2 + a'_4}, a'_5, a'_6)
$$
It follows that the map
$$
R(134)R(234): \BA_1 \lra \BA_3
$$
acts as
$$
(a''_1, \ldots, a''_6) = (a_1, 
\frac{a_3a_5a_6}{a_2a_4 + a_2a_6 + a_5a_6}, 
\frac{a_2a_4 + a_2a_6 + a_5a_6}{a_4 + a_6}, 
$$
$$
\frac{a_2a_3(a_4 + a_6)}{a_2a_4 + a_2a_6 + a_5a_6}, 
a_4 + a_6, \frac{a_4a_5}{a_4 + a_6}).
$$
The map
$$
L(1)L(4):\ \BA_3 \lra \BA_5
$$
interchanges $1$ with $2$ and $4$ with $5$, so the map
$$
C_{++} := L(1)L(4)R(134)R(234):\ \BA_1 \lra \BA_5
$$
acts as
$$
(a_1^{iv}, \ldots, a_6^{iv}) =  
(\frac{a_3a_5a_6}{a_2a_4 + a_2a_6 + a_5a_6}, a_1,  
\frac{a_2a_4 + a_2a_6 + a_5a_6}{a_4 + a_6}, 
$$
$$
a_4 + a_6, 
\frac{a_2a_3(a_4 + a_6)}{a_2a_4 + a_2a_6 + a_5a_6}, 
\frac{a_4a_5}{a_4 + a_6}).
$$

\subsection{$C_{+-}.$} The map
$$
R(124):\ \BA_5 \lra \BA_6
$$
acts as
$$
(a_1^v, \ldots, a_6^v) = (a_1^{iv}, \frac{a_3^{iv}a_4^{iv}}{a_2^{iv} + a_4^{iv}}, 
a_2^{iv} + a_4^{iv}, \frac{a_2^{iv}a_3^{iv}}{a_2^{iv} + a_4^{iv}}, 
a_5^{iv}, a_6^{iv}).
$$
The map
$$
R(123):\ \BA_6 \lra \BA_7
$$
acts as
$$
(a_1^{vi}, \ldots, a_6^{vi}) = (a_1^v, a_2^v, a_3^v, 
\frac{a_5^va_6^v}{a_4^v + a_6^v}, 
a_4^v + a_6^v, 
\frac{a_4^va_5^v}{a_4^v + a_6^v})
$$
It follows that the map
$$
R(123)R(124):\ \BA_5 \lra \BA_7
$$
acts as
$$
(a_1^{vi}, \ldots, a_6^{vi}) = (a_1^{iv}, 
\frac{a_3^{iv}a_4^{iv}}{a_2^{iv} + a_4^{iv}}, 
a_2^{iv} + a_4^{iv},
$$
$$
\frac{a_5^{iv}a_6^{iv}(a_2^{iv} + a_4^{iv})}
{a_2^{iv}a_3^{iv} + a_2^{iv}a_6^{iv} + a_4^{iv}a_6^{iv}},
\frac{a_2^{iv}a_3^{iv} + a_2^{iv}a_6^{iv} + a_4^{iv}a_6^{iv}}{a_2^{iv} + a_4^{iv}}, 
\frac{a_2^{iv}a_3^{iv}a_5^{iv}}{a_2^{iv}a_3^{iv} + a_2^{iv}a_6^{iv} + a_4^{iv}a_6^{iv}}).
$$
The map 
$$
L(3): \BA_7 \lra \BA_8
$$
interchanges $3$ with $4$, therefore
$$
C_{+-}:= L(3)R(123)R(124):\ \BA_5 \lra \BA_8
$$
acts as
$$
(a_1^{vii}, \ldots, a_6^{vii}) = (a_1^{iv}, 
\frac{a_3^{iv}a_4^{iv}}{a_2^{iv} + a_4^{iv}}, 
\frac{a_5^{iv}a_6^{iv}(a_2^{iv} + a_4^{iv})}
{a_2^{iv}a_3^{iv} + a_2^{iv}a_6^{iv} + a_4^{iv}a_6^{iv}},
$$
$$
a_2^{iv} + a_4^{iv}, 
\frac{a_2^{iv}a_3^{iv} + a_2^{iv}a_6^{iv} + a_4^{iv}a_6^{iv}}{a_2^{iv} + a_4^{iv}}, 
\frac{a_2^{iv}a_3^{iv}a_5^{iv}}{a_2^{iv}a_3^{iv} + a_2^{iv}a_6^{iv} + a_4^{iv}a_6^{iv}}).
$$ 

\subsection{LHS} The left hand side of the sonnet equation 
$$
C_+ := C_{+-}C_{++} = L(3)R(123)R(124)L(1)L(4)R(134)R(234): \BA_1 \lra \BA_8
$$
acts as follows
$$
= 
(\frac{a_3a_5a_6}{a_2a_4 + a_2a_6 + a_5a_6}, 
\frac{a_2a_4 + a_2a_6 + a_5a_6}{a_1 + a_4 + a_6}, 
\frac{a_2a_3a_4a_5(a_1 + a_4 + a_6)}
{(a_1a_2 + a_1a_5 + a_4a_5)(a_2a_4 + a_2a_6 + a_5a_6)},
$$
$$ 
a_1 + a_4 + a_6,
\frac{a_1a_2 + a_1a_5 + a_4a_5}
{a_1 + a_4 + a_6}, 
\frac{a_1a_2a_3}{a_1a_2 + a_1a_5 + a_4a_5})
$$ 
We have used two identites
$$
a_1a_2a_4 + a_1a_2a_6 +a_1a_5a_6 +a_1a_4a_5 +a_4a_4a_5 +a_4a_5a_6 = 
$$
\begin{equation}
= 
(a_4 + a_6)(a_1a_2 + a_1a_5 + a_4a_5),  
\end{equation}
and
$$
a_1a_2a_4 + a_2a_4a_4 + a_2a_4a_6 + a_1a_2a_6 + a_1a_5a_6 + a_4a_5a_6 =  
$$
\begin{equation}
= (a_1 + a_4)(a_2a_4 + a_2a_6 + a_5a_6) 
\end{equation}
The first one implies that 
$$
a_2^{iv}a_3^{iv} + a_2^{iv}a_6^{iv} + a_4^{iv}a_6^{iv}
= a_1a_2 + a_1a_5 + a_4a_5
$$

\

It folows that 
$$
(a_1^{vii}, \ldots, a_6^{vii}) = 
$$
$$
= 
(\frac{a_3a_5a_6}{a_2a_4 + a_2a_6 + a_5a_6}, 
\frac{a_2a_4 + a_2a_6 + a_5a_6}{a_1 + a_4 + a_6}, 
\frac{a_2a_3a_4a_5(a_1 + a_4 + a_6)}
{(a_1a_2 + a_1a_5 + a_4a_5)(a_2a_4 + a_2a_6 + a_5a_6)},
$$
$$ 
a_1 + a_4 + a_6,
\frac{a_1a_2 + a_1a_5 + a_4a_5}
{a_1 + a_4 + a_6}, 
\frac{a_1a_2a_3}{a_1a_2 + a_1a_5 + a_4a_5})
$$

\subsection{The lower road. $C_{-+}$.} The map 
$$
R(123):\ \BB_1 \lra \BB_2
$$
acts as follows
$$
(a_1'', \ldots, a_6'') = (\frac{a_2'a_3'}{a_1' + a_3'}, a_1' + a_3', 
\frac{a_1'a_2'}{a_1' + a_3'}, 
a_4', a_5', a_6') 
$$
The map  
$$
R(124):\ \BB_2 \lra \BB_3
$$
acts as follows
$$
(a_1''', \ldots, a_6''') = (a_1'', a_2'', 
\frac{a_4''a_5''}{a_3'' + a_5''}, a_3'' + a_5'', 
\frac{a_3''a_4''}{a_3'' + a_5''}, a_6'') 
$$
Therefore the map
$$
R(124)R(123):\ \BB_1 \lra \BB_3
$$
acts as follows
$$
(a_1''', \ldots, a_6''') = (\frac{a_2'a_3'}{a_1' + a_3'}, 
a_1' + a_3', \frac{a_4'a_5'(a_1' + a_3')}{a_1'a_2' + a_1'a_5' + a_3'a_5'}, 
$$
$$
\frac{a_1'a_2' + a_1'a_5' + a_3'a_5'}{a_1' + a_3'},
\frac{a_1'a_2'a_4'}{a_1'a_2' + a_1'a_5' + a_3'a_5'}, a_6')
$$
The map 
$$
L(3):\ \BB_0 \lra \BB_1 
$$
interchanges $3$ with $4$. It follows that the map
$$
C_{-+}:= R(124)R(123)L(3):\ \BB_0 \lra \BB_3
$$
acts as 
$$
(a_1''', \ldots, a_6''') = (\frac{a_2a_4}{a_1 + a_4}, 
a_1 + a_4, \frac{a_3a_5(a_1 + a_4)}{a_1a_2 + a_1a_5 + a_4a_5}, 
$$
$$
\frac{a_1a_2 + a_1a_5 + a_4a_5}{a_1 + a_4},
\frac{a_1a_2a_3}{a_1a_2 + a_1a_5 + a_4a_5}, a_6)
$$

\subsection{$C_{--}$.} The map
$$
R(134): \BB_5 \lra \BB_6
$$
acts as follows
$$
(a_1^{vi}, \ldots, a_6^{vi}) = (a_1^v, a_2^v, \frac{a_4^va_5^v}{a_3^v + a_5^v}, 
a_3^v + a_5^v, 
\frac{a_3^va_4^v}{a_3^v + a_5^v}, a_6^v) 
$$
The map
$$
R(234): \BB_6 \lra \BB_7
$$
acts as follows
$$
(a_1^{vii}, \ldots, a_6^{vii}) = (\frac{a_2^{vi}a_3^{vi}}{a_1^{vi} + a_3^{vi}}, 
a_1^{vi} + a_3^{vi}, \frac{a_1^{vi}a_2^{vi}}{a_1^{vi} + a_3^{vi}},
a_4^{vi}, a_5^{vi}, a_6^{vi}).   
$$
It follows that the map
$$
R(234)R(134): \BB_5 \lra \BB_7
$$
acts as follows
$$
(a_1^{vii}, \ldots, a_6^{vii}) = (
\frac{a_2^va_4^va_5^v}{a_1^va_3^v + a_1^va_5^v + a_4^va_5^v}, 
\frac{a_1^va_3^v + a_1^va_5^v + a_4^va_5^v}{a_3^v + a_5^v}, 
\frac{a_1^va_2^v(a_3^v + a_5^v)}{a_1^va_3^v + a_1^va_5^v + a_4^va_5^v}, 
$$
$$
a_3^v + a_5^v, \frac{a_3^va_4^v}{a_3^v + a_5^v}, a_6^v).
$$
The map 
$$
L(5)L(2): \BB_3 \lra \BB_5
$$
interchanges $2$ with $3$ and $5$ with $6$. 

It follows that 
$$
C_{--} := R(234)R(134)L(5)L(2): \BB_3 \lra \BB_7
$$
acts as 
$$
(a_1^{vii}, \ldots, a_6^{vii}) = (
\frac{a_3'''a_4'''a_6'''}{a_1'''a_2''' + a_1'''a_6''' + a_4'''a_6'''}, 
\frac{a_1'''a_2''' + a_1'''a_6''' + a_4'''a_6'''}{a_2''' + a_6'''}, 
\frac{a_1'''a_3'''(a_2''' + a_6''')}{a_1'''a_2''' + a_1'''a_6''' + a_4'''a_6'''}, 
$$
$$
a_2''' + a_6''', \frac{a_2'''a_4'''}{a_2''' + a_6'''}, a_5''').
$$

\subsection{RHS.} The right hand side of the sonnet equation
$$
C_- := C_{--}C_{-+} = R(234)R(134)L(5)L(2)R(124)R(123)L(3): \BB_0 \lra \BB_7
$$
acts as 
$$
(a_1^{vii}, \ldots, a_6^{vii}) = 
$$
$$
= (\frac{a_3a_5a_6}{a_2a_4 + a_2a_6 + a_5a_6},
\frac{a_2a_4 + a_2a_6 + a_5a_6}{a_1 + a_4 + a_6}, 
\frac{a_2a_3a_4a_5(a_1 + a_4 + a_6)}{(a_1a_2 + a_1a_5 + a_4a_6)(a_2a_4 + a_2a_6 + a_5a_6)},
$$
$$
a_1 + a_4 + a_6, \frac{a_1a_2 + a_1a_5 + a_4a_6}{a_1 + a_4 + a_6}, 
\frac{a_1a_2a_3}{a_1a_2 + a_1a_5 + a_4a_6})
$$

We have used the identity
$$
a_1a_2a_4 + a_2a_4a_4 + a_2a_4a_6 + a_1a_2a_6 + a_1a_5a_6 + a_4a_5a_6 =  
$$
\begin{equation}
= (a_1 + a_4)(a_2a_4 + a_2a_6 + a_5a_6) 
\end{equation}
which implies that
\begin{equation}
a_1'''a_2''' + a_1'''a_6''' + a_4'''a_6'''
= a_2a_4 + a_2a_6 + a_5a_6
\end{equation}

We see that $C_+ = C_-$, i.e. the sonnet equation holds true.








\section{Two parametric upper diagonal case}\label{flacon}

In this Section we propose an example which encompasses the Lusztig and the Berenstein - Zelevinsky flips. 

\subsection{} Consider a matrix $\begin{pmatrix}a^{-1}&x\cr 0 &a\end{pmatrix}$.

For the decomposition $121$ we have 
\[
\begin{pmatrix}
a^{-1}&x&0 \cr 0 &a&0\cr
0&0&1
\end{pmatrix}
\begin{pmatrix}
 1  &0&0\cr 
0& b^{-1}&y\cr
0&0&b
\end{pmatrix}
\begin{pmatrix}
c^{-1}&z&0 \cr 0 &c&0\cr
0&0&1
\end{pmatrix}=\]
\begin{equation}\label{delta1}
 \begin{pmatrix}
\frac1{ac}&\frac{cx}{b}+\frac{z}{a}&xy \cr 0 &\frac {ac}{b}&ay\cr
0&0&b
\end{pmatrix}
\end{equation}

For the decomposition $212$ we have

\[
\begin{pmatrix}
 1  &0&0\cr 
0& a'^{-1}&x'\cr
0&0&a'
\end{pmatrix}
\begin{pmatrix}
b'^{-1}&y'&0 \cr 0 &b'&0\cr
0&0&1
\end{pmatrix}
\begin{pmatrix}
 1  &0&0\cr 
0& c'^{-1}&z'\cr
0&0&c'
\end{pmatrix}
=\]
\begin{equation}\label{delta2}
 \begin{pmatrix}
\frac 1{b'}&\frac{y'}{c'} & y'z' \\
0 & \frac{b'}{a'c'}&\frac{b'z'}{a'}+c'x'\cr
0&0&a'c'
\end{pmatrix}
\end{equation}

Hence the transformation $121\to 212$ gives us  the correspondence 
$$
(a, b, c, x,y,z)\to 
$$
\begin{equation}\label{bzL1}
\to (b/c', \, ac, \, c', \,\frac{zyb}{c'}\frac 1{za^{-1}b+xc}, \,\frac{c'}{b} (za^{-1}b+xc) , \,\frac{xyb}{c'}\frac 1{za^{-1}b+xc}).
\end{equation}

Note that if in (\ref{bzL1}) we set $a=b=c=c'=1$, we get the Lusztig transformation 
\[
(1,1,1,x,y,z)\to (1,1,1,\frac{zy}{x+z}, x+z, \frac{xy}{x+z}).
\]
If we put $x=y=z=1$ and $c'=a+\frac{b}{c}$, we get the Sergeev's case $(\beta)$ or BZ-transformation for upper triangular version, that  is one has to swap $a$ and $c$ in the BZ-transformation,
\[
(a,b,c,1,1,1)\to (\frac{bc}{a+bc}, ac, a+\frac{b}{c}, 1, 1, 1).
\]

The correspondence (\ref{bzL1})  partially factorizes,  that is one can consider a correspondence
\begin{equation}\label{BZnew}
(a,b,c)\to (\frac {b}{c'}, ac, c').
\end{equation}

The upper path for (\ref{BZnew})  is following 
\[
(a,b,c,d,e,f)\to (\frac{b}{c'}, ac, c',d,e,f)\to (\frac{b}{c'}, ac, \frac{d}{e'},c'e,e',f)\to (\frac{b}{c'}, \frac{d}{e'}, acf, \frac{c'e}{f'},f',e')\to 
\]
\begin{equation}\label{upper}
(\frac{d}{d'e'},  \frac{be}{f'}, acf, \overline d, f', e')
\end{equation}

The lower path  is
\[
(a,b,c,d,e,f)\to (a,b, d,\frac{e}{f_1},cf,f_1)\to (a,\frac {d}{e_1},\frac{be}{f_1}, {e_1},cf,f_1)\to (\frac {d}{e'_1}, \frac{be}{f_1d_1},acf, {d_1},e_1, f_1)\to
\]
\begin{equation}\label{lower}
(\frac {d}{e_1}, \frac{be}{f_1d_1}, acf, \frac{e_1}{{\overline f_1}}, {d_1f_1}, \overline f_1)
\end{equation}

Then (\ref{upper}) $=$ (\ref{lower}) reads as follows 

\[
\frac{d}{e_1}= \frac{d}{d'e'},\]
\[
\frac{be}{f_1d_1}=\frac{be}{f'},
\]
\[
acf=acf,
\]
\[
\frac {e_1}{\overline f_1}=\overline d,
\]
\[
c_1f_1=f' ,
\]
\[
\overline f_1=e' .
\]

From this system we 
get
\begin{equation}\label{uppersep1}
f'=c_1f_1,\, e'=\overline f_1, \, \overline d=\frac {e_1}{\overline f_1}.
\end{equation}

This system is obtained as the equality of the upper and the lower paths for $(a,b,c,d,e,f)$ coordinates.

\subsection{The upper path 2}

\

Here we compute the upper path for the transformation of
\[
(x,y,z,u,v,w)
\]
under (\ref{bzL1}).
The resulting  6-tuple takes the form 
\begin{equation}\label{upper2}
(\frac 1{e'\overline d}P_1, \frac{\overline de'}{f'}P_2, \frac {f'}{e'^2\overline d}P_3, f'P_4, , \frac{e'}{f'}P_5, \frac 1{e'}P_6), 
\end{equation}
where $P_i$ are rational functions in $(a,b,c,d,e,f,x,y,z,u,v,w)$.

\subsection{The lower path 2}

\

The lower path sends $(x,y,z,u,v,w)$ to the 6-tuple 
\begin{equation}\label{upper2}
(\frac 1{e_1}R_1, \frac{e_1}{f_1d_1}R_2, \frac {f_1d_1}{e_1\overline f_1}R_3, {f_1d_1}R_4, \frac {\overline f_1}{d_1f_1}R_5, \frac{1}{\overline f_1}R_6), 
\end{equation}
where $R_i$ are rational functions in $(a,b,c,d,e,f,x,y,z,u,v,w)$.

For the consistenccy of these paths, taking into account (\ref{uppersep1}), we get the following  six equations for $12$ variables
$(a, b, c,d,e,f, x,y,z,uv,w)$,
\[
P_i=R_i, \qquad i=1, \cdots , 6.
\]


\


\section{Tetrahedron equations for some particular  cases}

(a) {\em Lusztig flip}

\subsection{•} The aim of this Section is to provide a proof of \cite{SV}, Theorem 8.1. 
It says that the diagram (\ref{tetr-eq-diagram}) commutes, where $R$ stands for the Lusztig map 
$$
(a,b,c)\mapsto (bc/(a+c), a+c, ab/(a+c)),
$$
cf. \cite{L}, Proof of Prop. 42.2.4, p. 336.

\

We use the notations from \cite{SV}, 8.2.2.

\subsection{• Notations and identities.} Linear:
$$
\alpha = a_1 + a_3; \gamma = a_3 + a_6;\ \beta = a_1 +  a_3 + a_6;
$$
Quadratic:
$$
\delta = a_2a_3 + a_2a_6 + a_5a_6;\  
\epsilon = a_1a_2 + a_1a_5 + a_3a_5.
$$

Cubical identities:
\begin{equation}\label{cubic1}
\mu:= a_6\epsilon + a_2a_3\beta = \alpha\delta,
\end{equation}
\begin{equation}\label{cubic2}
\nu: = a_1\delta + a_3a_5\beta = \gamma\epsilon.
\end{equation}

GO!

\subsection{•}
On the one hand: the way
$$
\begin{matrix}\cdot & \lra & \cdot\\ \uparrow & & \downarrow 
\end{matrix}
$$
gives
$$
\left(\begin{matrix} 1 & 2 & 3 & 1 & 2 & 1\\ 
a_6 & a_5 & a_4 & a_3 & a_2 & a_1 \end{matrix}\right)\ \overset{L(3)}\lra
$$
$$
\left(\begin{matrix} 1 & 2 & 1 & 3 & 2 & 1\\ 
a_6 & a_5 & a_3 & a_4 & a_2 & a_1 \end{matrix}\right)\ \overset{R(1)}\lra
$$ 
$$
\left(\begin{matrix} 2 & 1 & 2 & 3 & 2 & 1\\ 
a_3a_5/\gamma & \gamma & a_5a_6/\gamma & a_4 & a_2 & a_1 \end{matrix}\right)\ \overset{R(3)}\lra
$$
$$
\left(\begin{matrix} 2 & 1 & 3 & 2 & 3 & 1\\ 
a_3a_5/\gamma & \gamma & a_2a_4\gamma/\delta & \delta/\gamma & a_4a_5a_6/\delta  & a_1 \end{matrix}\right)\ \overset{L(2)L(5)}\lra
$$
$$
\left(\begin{matrix} 2 & 3 & 1 & 2 & 1 & 3\\ 
a_3a_5/\gamma & a_2a_4\gamma/\delta & \gamma & \delta/\gamma & a_1 & a_4a_5a_6/\delta \end{matrix}\right)\ \overset{R(3)}\lra
$$
$$
\left(\begin{matrix} 2 & 3 & 2 & 1 & 2 & 3\\ 
a_3a_5/\gamma & a_2a_4\gamma/\delta & a_1\delta/\gamma\beta & \beta & \delta/\beta & a_4a_5a_6/\delta \end{matrix}\right)\ \overset{R(1)}\lra
$$
$$
\left(\begin{matrix} 3 & 2 & 3 & 1 & 2 & 3\\ 
a_1a_2a_4/\epsilon & \epsilon/\beta & a_2a_3a_4a_5\beta/\delta\epsilon & \beta & \delta/\beta & a_4a_5a_6/\delta \end{matrix}\right)\ 
$$
On the last step one uses the cubical identity (\ref{cubic2}).

\subsection{•} On the other hand, the way
$$
\begin{matrix}\downarrow & & \uparrow\\ \cdot & \lra & \cdot 
\end{matrix}
$$
gives
$$
\left(\begin{matrix} 1 & 2 & 3 & 1 & 2 & 1\\ 
a_6 & a_5 & a_4 & a_3 & a_2 & a_1 \end{matrix}\right)\ \overset{R(4)}\lra
$$
$$
\left(\begin{matrix} 1 & 2 & 3 & 2 & 1 & 2\\ 
a_6 & a_5 & a_4 & a_1a_2/\alpha & \alpha & a_2a_3/\alpha \end{matrix}\right)\ 
\overset{R(2)}\lra
$$
$$
\left(\begin{matrix} 1 & 3 & 2 & 3 & 1 & 2\\ 
a_6 & a_1a_2a_4/\epsilon & \epsilon/\alpha & a_4a_5\alpha/\epsilon & \alpha & a_2a_3/\alpha \end{matrix}\right)\ 
\overset{L(1)L(4)}\lra
$$
$$
\left(\begin{matrix} 3 & 1 & 2 & 1 & 3 & 2\\ 
a_1a_2a_4/\epsilon & a_6 & \epsilon/\alpha & \alpha & a_4a_5\alpha/\epsilon & a_2a_3/\alpha \end{matrix}\right)\ 
\overset{R(2)}\lra
$$
$$
\left(\begin{matrix} 3 & 2 & 1 & 2 & 3 & 2\\ 
 a_1a_2a_4/\epsilon & \epsilon/\beta & \beta & a_6\epsilon/\alpha\beta & a_4a_5\alpha/\epsilon & a_2a_3/\alpha \end{matrix}\right)\ 
\overset{R(4)}\lra
$$
$$
\left(\begin{matrix} 3 & 2 & 1 & 3 & 2 & 3\\ 
 a_1a_2a_4/\epsilon & \epsilon/\beta & \beta & a_2a_3a_4a_5\beta/\epsilon\delta & 
 \delta/\beta  & a_4a_5a_6/\delta \end{matrix}\right)\ 
 \overset{L(3)}\lra
$$
$$
\left(\begin{matrix} 3 & 2 & 3 & 1 & 2 & 3\\ 
a_1a_2a_4/\epsilon & \epsilon/\beta & a_2a_3a_4a_5\beta/\delta\epsilon & \beta & \delta/\beta & a_4a_5a_6/\delta \end{matrix}\right)\ 
$$
On the last step one uses the cubical identity (\ref{cubic1}).

\

Two ways have given the same answer, which proves our assertion. $\Box$

\

For another proof, see \cite{L}, Prop. 42.2.4. 

\

(b) {\em Sergeev $(\alpha)$ flip}

\subsection{•} Consider another involution
$$
\bR:\ \BC(a, b, c) \lra \BC(a, b, c)
$$
defined by
$$
a' = c,\ b' = ac - b,\ c' = a.
$$
It coincides with \cite{S}, (17) with $k = - 1$.

\subsection{Useful quantities.}
$$
\alpha = \delta' = a_1a_3 - a_2;\ \alpha' = \delta = a_3a_6 - a_5
$$
$$
\beta = a_1a_5 - a_4;\ \beta' = a_2a_6 - a_4
$$
$$
\gamma = \gamma' = a_6\alpha - \beta = a_1\alpha' - \beta' = 
$$
$$
 = a_1a_3a_6 - a_1a_5 - a_2a_6 + a_4.
$$
A cubical relation:
$$
a_6 \alpha - \beta = a'_1 \alpha' - \beta'
\eqno{(2.2.1)}
$$

\subsection{Lemma.} 
$$
\gamma = \left|\begin{matrix}
a_1 & a_2 & a_4\\
1 & a_3 & a_5\\
0 & 1 & a_6
\end{matrix}\right|.
$$
$\Box$

\subsection{ Evolution 1.} The way
$$
\begin{matrix}\cdot & \lra & \cdot\\ \uparrow & & \downarrow 
\end{matrix}
$$
gives
$$
\left(\begin{matrix} 1 & 2 & 3 & 1 & 2 & 1\\ 
a_6 & a_5 & a_4 & a_3 & a_2 & a_1 \end{matrix}\right)\ \overset{L(3)}\lra
$$
$$
\left(\begin{matrix} 1 & 2 & 1 & 3 & 2 & 1\\ 
a_6 & a_5 & a_3 & a_4 & a_2 & a_1 \end{matrix}\right)\ \overset{R(1)}\lra
$$ 
$$
\left(\begin{matrix} 2 & 1 & 2 & 3 & 2 & 1\\ 
a_3 & \alpha' & a_6 & a_4  & a_2 & a_1 \end{matrix}\right)\ \overset{R(3)}\lra
$$
$$
\left(\begin{matrix} 2 & 1 & 3 & 2 & 3 & 1\\ 
a_3 & \alpha' & a_2 & \beta' & a_6  & a_1 \end{matrix}\right)\ \overset{L(2)L(5)}\lra
$$
$$
\left(\begin{matrix} 2 & 3 & 1 & 2 & 1 & 3\\ 
a_3 & a_2 & \alpha' & \beta' & a_1 & a_6 
\end{matrix}\right)\ \overset{R(3)}\lra
$$
$$
\left(\begin{matrix} 2 & 3 & 2 & 1 & 2 & 3\\ 
a_3 & a_2 & a_1 & \gamma' & \alpha' & a_6 
\end{matrix}\right)\ \overset{R(1)}\lra
$$
$$
\left(\begin{matrix} 3 & 2 & 3 & 1 & 2 & 3\\ 
a_1 & \delta' & a_3 & \gamma' & \alpha' & a_6
\end{matrix}\right)\ 
$$

\subsection{ Evolution 2.} The way
$$
\begin{matrix}\downarrow & & \uparrow\\ \cdot & \lra & \cdot 
\end{matrix}
$$
gives
$$
\left(\begin{matrix} 1 & 2 & 3 & 1 & 2 & 1\\ 
a_6 & a_5 & a_4 & a_3 & a_2 & a_1 \end{matrix}\right)\ \overset{R(4)}\lra
$$
$$
\left(\begin{matrix} 1 & 2 & 3 & 2 & 1 & 2\\ 
a_6 & a_5 & a_4 & a_1 & \alpha & a_3 \end{matrix}\right)\ 
\overset{R(2)}\lra
$$
$$
\left(\begin{matrix} 1 & 3 & 2 & 3 & 1 & 2\\ 
a_6 & a_1 & \beta & a_5 & \alpha & a_3 
\end{matrix}\right)\ 
\overset{L(1)L(4)}\lra
$$
$$
\left(\begin{matrix} 3 & 1 & 2 & 1 & 3 & 2\\ 
a_1 & a_6 & \beta  & \alpha & a_5 & a_3
\end{matrix}\right)\ 
\overset{R(2)}\lra
$$
$$
\left(\begin{matrix} 3 & 2 & 1 & 2 & 3 & 2\\ 
a_1 & \alpha & \gamma & a_6 & a_5 & a_3 
\end{matrix}\right)\ 
\overset{R(4)}\lra
$$
$$
\left(\begin{matrix} 3 & 2 & 1 & 3 & 2 & 3\\ 
a_1 & \alpha & \gamma & a_3 & \delta & a_6 
\end{matrix}\right)\ 
\overset{L(3)}\lra
$$
$$
\left(\begin{matrix} 3 & 2 & 3 & 1 & 2 & 3\\ 
a_1 & \alpha & a_3 & \gamma  & \delta & a_6
\end{matrix}\right)\ 
$$
Now we recall that $\alpha = \delta', \delta = \alpha'$, and $\gamma = \gamma'$, whence both ways arrive at the same result.  

\subsection{•} Therefore it is natural to expect the in general in the evolution 
appear the minors of the unversal triangular matrix.

\section{Long products and $S(n)$ evolutions.}

\subsection{$c(n)$-fold products.} Pick $n$, so we have matrices 
$A(a,b,c)_i\in Mat_n$, \newline $1\leq n - 1$.  

\

Let us pick Coxeter generators $s_1, \ldots, s_{n-1}$ 
of $S(n)$ and consider a reduced decomposition of $w_0(n)$ into a product of $s_i$
\begin{equation}\label{reduced-dec}
w_0(n) = s_{i_1}s_{i_2}\ldots s_{i_{c(n)}},\ c(n) = n(n-1)/2. 
\end{equation}
Given $3c(n)$ numbers
\begin{equation}\label{numbers-abc}
\mathbf{a} = (a_1, b_1, c_1, \ldots, a_{c(n)}, b_{c(n)}, c_{c(n)})
\end{equation}
consider the product
$$
A(\mathbf{a}) = A(a_1,b_1,c_1)_{i_1}A(a_2,b_2,c_2)_{i_2}
\ldots A(a_{c(n)},b_{c(n)},c_{c(n)})_{i_{c(n)}}\in Mat_n
$$

\subsection{Theorem.} {\em For any reduced decomposition (\ref{reduced-dec}) this 
$(n\times n)$-matrix is $c$-upper triangular, $A(\mathbf{a})\in B^c_+$.}

\subsection{Truncated products.} For any $k$, $1\leq k\leq c(n)$, denote
$$
w_k = w_0(n)_{\leq k} = s_{i_1}\ldots s_{i_k},
$$
and
$$
w'_k = w_0(n)_{\geq k+1} = s_{i_{k+1}}\ldots s_{i_{c(n)}}
$$
Let
$$
A(\mathbf{a})_{\leq k} = A(a_1,b_1,c_1)_{i_1}A(a_2,b_2,c_2)_{i_2}
\ldots A(a_k,b_k,c_k)_{i_k}. 
$$  
Among reduced decompositions (\ref{reduced-dec}) there are two distinguished ones: the minimal one of the form
$$
w_0(n) = s_1s_2s_1s_3s_2s_1\ldots,
$$
and the maximal one having the form
$$
w_0(n) = s_ns_{n-1}s_ns_{n-2}s_{n-1}s_n\ldots .
$$

\ 

The following assertion generalizes (\ref{double-one}).
 
\subsection{Theorem.}\label{permut-decomp}{\em For any $1\leq k\leq c(n)$ consider the matrix $A(\mathbf{a})_{\leq k}$ corresponding to the minimal or to the maximal reduced decomposition of $w_0(n)$.} 

(a) {\em $w_k$ is a unique 
permutation $w_k\in S(n)$ such that $\iota(w_k)A(\mathbf{a})_{\leq k}\in B_+^c(n)$.}

(b) {\em $w'_k$ is a unique 
permutation $w'_k\in S(n)$ such that $A(\mathbf{a})_{\leq k}\iota(w'_k)\in B_+^c(n)$.} 

\

{\bf Idea of proof.} It is sufficient to suppose that all $a_i = b_i = c_i = 1$.

\

Moreover the following more general statement holds true:

\subsection{Theorem.}\label{permut-decomp-gen}{\em Let us pick an arbitrary reduced decomposition of $w_0(n)$. For any $1\leq k\leq n$   
consider the corresponding  matrix $A(\mathbf{a})_{\leq k}$.} 

(a) {\em  $w_k$ is a unique 
permutation $w_k\in S(n)$ such that $\iota(w_k)A(\mathbf{a})_{\leq k}\in B_+^c(n)$.}

(b) {\em $w'_k$ is a unique 
permutation $w'_k\in S(n)$ such that $A(\mathbf{a})_{\leq k}\iota(w'_k)\in B_+^c(n)$.}

\section{Wronskian evolutions}\label{sec-wronskian-ev}

In this Section we suppose that our ground ring contains $\BQ$. 

\subsection{Start of procedure.} The above considerations allow to describe the corresponding Wronskian evolutions. There are four of them: upper- and lower triangular, and $c$-upper and $c$-lower triangular evolutions. 

\

Let us consider the $c$-upper triangular case. Fix an ineger $n\geq 1$. Similarly to \ref{triple-cont} and \ref{quadruple} 
we start with a collection of polynomials 
$u_i(x) \in \BC[x]$, $u_i(x) = x^{i-1}/(i-1)! + \ldots$, $1\leq i \leq n$. Here 
$\ldots$ means the terms of higher degree. 

\

Fix a reduced decomposition (\ref{reduced-dec}). Pick numbers (\ref{numbers-abc}). We have
$$
A(\mathbf{a})
\left(\begin{matrix}
u_n \\ u_{n-1} \\ . \\ u_1 
\end{matrix}\right) = \left(\begin{matrix}
d_1 + \ldots \\ d_2 x + \ldots \\ . \\ d_n x^{n-1}/(n-1)! + \ldots 
\end{matrix}\right),\ d_i\in \BC^*,
$$
cf. \ref{quadruple-eq}. 

\subsection{Wronskian coordinates.} For any $1\leq k \leq n$ 
 we have, according to Theorem \ref{permut-decomp-gen}, the element $w_k\in S(n)$ such that $\iota(w_k)A(\mathbf{a})_{\leq k}\in B_+^c(n)$. 

Consider the column vector
$$
\left(\begin{matrix}
v_{1k} \\ v_{2k} \\ . \\ v_{nk} 
\end{matrix}\right) = \iota(w_k)A(\mathbf{a})_{\leq k}\left(\begin{matrix}
u_n \\ u_{n-1} \\ . \\ u_1 
\end{matrix}\right)
$$

\


We define
$$
w_k(x) := Wr(v_{1k}(x), \ldots, v_{kk}(x)). 
$$
Then we have 
$$
w_k(x) = e_i + \ldots, \ e_i\in \BC^*.
$$

\

G.K.: Kharkevich Institute for Information Transmission Problems, Moscow,
\newline Russia 127994; 
{koshevoyga@gmail.com}

V.S.: Institut de Math\'ematiques de Toulouse, 
Universit\'e Toulouse III Paul Sabatier, 119 Route de Narbonne, 31062 Toulouse, France;
Kavli IPMU, 5-1-5 Kashiwanoha, Kashiwa, Chiba, 277-8583 Japan; 
{schechtman@math.ups-tlse.fr}

A.V.: Department of Mathematics, University of North Carolina at Chapel Hill,
Chapel Hill, NC 27599-3250, USA;\ Faculty of Mathematics and Mechanics, Lomonosov
Moscow State University, Leninskiye Gory 1, 119991 Moscow GSP-1, 
\newline Russia; 
{anv@email.unc.edu}

\end{document}